\documentclass[leqno,12pt]{article}  
\usepackage{amsmath}
\usepackage{amssymb}

\usepackage[german,english]{babel}
\usepackage{amsthm}

\title{Integral structures in automorphic line bundles on the $p$-adic upper half plane}

\author{\textsc{Elmar Grosse-Kl\"onne}}
\date{}

\theoremstyle{plain} 
\newtheorem{satz}{Theorem}[section]  
\newtheorem{kor}[satz]{Corollary}  
\newtheorem{lem}[satz]{Lemma}  
\newtheorem{pro}[satz]{Proposition}  
\newcommand{\ho}{\mbox{\rm Hom}}  

\newcommand{\spec}{\mbox{\rm Spec}}  

\newcommand{\spf}{\mbox{\rm Spf}}  

\newcommand{\bi}{\mbox{\rm Im}}  
\newcommand{\ke}{\mbox{\rm Ker}}  

\newcommand{\kara}{\mbox{\rm char}}  

\newcommand{\sym}{\mbox{\rm Sym}}
\newcommand{\id}{\mbox{\rm id}}

\newcommand{\dlog}{\mbox{\rm dlog}}

\newcommand{\ind}{\mbox{\rm Ind}}

\theoremstyle{remark}

\theoremstyle{definition}

\DeclareMathOperator{\Hom}{Hom}

\newcommand{\0}{\ensuremath{\overrightarrow{0}}}

\begin{document}
\maketitle
\footnote[0]{I wish to thank Peter Schneider and Jeremy Teitelbaum for generously providing me with some helpful private notes on their own work, and for their interest. I am also grateful to Matthias Strauch for useful discussions on odd weight modular forms. I thank Christophe Breuil for his interest and his insisting on lattices for the entire $G$-action. Finally I thank the referee for his suggestions concerning the presentaton of several technical constructions.}

\begin{abstract}
Given an automorphic line bundle ${\mathcal O}_X(k)$ of weight $k$ on the Drinfel'd upper half plane $X$ over a local field $K$, we construct a ${\rm GL}\sb 2(K)$-equivariant integral lattice ${\mathcal O}_{\widehat{\mathfrak X}}(k)$ in ${\mathcal O}_X(k)\otimes_K\widehat{K}$, as a coherent sheaf on the formal model $\widehat{\mathfrak{X}}$ underlying $X\otimes_K\widehat{K}$. Here $\widehat{K}/K$ is ramified of degree $2$. This generalizes  a construction of Teitelbaum from the case of even weight $k$ to arbitrary integer weight $k$. We compute $H^*(\widetilde{\mathfrak{X}},{\mathcal O}_{\widehat{\mathfrak X}}(k))$ and obtain applications to the de Rham cohomology $H_{dR}^1(\Gamma\backslash X,\sym_K^k({\rm St}))$ with coefficients in the $k$-th symmetric power of the standard representation of ${\rm SL}\sb 2(K)$ (where $k\ge0$) of projective curves $\Gamma\backslash X$ uniformized by $X$: namely, we prove the degeneration of a certain reduced Hodge spectral sequence computing $H_{dR}^1(\Gamma\backslash X,\sym_K^k({\rm St}))$, we re-prove the Hodge decomposition of $H_{dR}^1(\Gamma\backslash X,\sym_K^k({\rm St}))$ and show that the monodromy operator on $H_{dR}^1(\Gamma\backslash X,\sym_K^k({\rm St}))$ respects integral de Rham structures and is induced by a "universal"{} monodromy operator defined on $\widehat{\mathfrak{X}}$, i.e. before passing to the $\Gamma$-quotient.  
\end{abstract}

\begin{center} {\bf Introduction}
\end{center} Let $K$ be a local field and let $X$ be the Drinfel'd upper half plane over $K$; that is, the projective line over $K$ with its $K$-rational points removed. $G={\rm GL}\sb 2(K)$ acts on $X$. Let ${\mathcal O}_X(k)$ be the structure sheaf on the rigid space $X$, endowed with the automorphic action by $G$ of weight $k\in\mathbb{Z}$. For $k\ge0$ and even, Teitelbaum \cite{jer} constructed a $G$-invariant integral lattice in ${\mathcal O}_X(k)$, as a line bundle on the natural formal ${\mathcal O}_K$-scheme ${\mathfrak X}$ underlying $X$. He then reduced this bundle modulo the maximal ideal of ${\mathcal O}_K$ and determined explicitly its global sections, as a representation of $G$ on an infinite dimensional vector space over the residue field  ${\mathbb{F}}$ of $K$. The first aim of this paper is to extend his results to {\it any} weight $k\in\mathbb{Z}$. Now it is not hard to see that for odd $k$ there is no $G$-equivariant ${\mathcal O}_{\mathfrak{X}}$-line bundle lattice in ${\mathcal O}_X(k)$. Let $\widehat{K}$ be a ramified extension of $K$ of degree $2$, let $\widehat{\mathfrak X}={\mathfrak X}\otimes_{{\mathcal O}_K}{\mathcal O}_{\widehat{K}}$ be the base extended formal ${\mathcal O}_{\widehat{K}}$-scheme, let $\widetilde{\mathfrak X}={\mathfrak X}\otimes_{{\mathcal O}_K}{\mathbb{F}}=\widehat{\mathfrak X}\otimes_{{\mathcal O}_{\widehat{K}}}{\mathbb{F}}$. We show that for any $k\in\mathbb{Z}$, if we twist the automorphic action on ${\mathcal O}_X(k)$ by a suitable character, there is a $G$-equivariant ${\cal O}_{\widehat{\mathfrak X}}$-module ${\mathcal O}_{\widehat{\mathfrak X}}(k)$ which is a lattice inside ${\mathcal O}_X(k)\otimes_K\widehat{K}$. If $k$ is even it is a line bundle, if $k$ is odd it is not: around the singular points of $\widehat{\mathfrak X}$ it needs two generators. We show that $H^0(\widetilde{\mathfrak X},{\mathcal O}_{\widehat{\mathfrak X}}(k))$ for $k\ge0, k\ne 1$ and $H^1(\widetilde{\mathfrak X},{\mathcal O}_{\widehat{\mathfrak X}}(k))$ for $k\le-1$ are {\it precisely} those cohomology groups which do not vanish. We prove that they are ${\mathcal O}_{\widehat{K}}$-flat and that their formation commutes with base change to the special fibre  $\widetilde{\mathfrak X}$. We determine these $G$-representations obtained by reduction modulo the maximal ideal of ${\mathcal O}_{\widehat{K}}$ in the same manner as in \cite{jer}. Next we establish for $k\ge2$ an isomorphism between $H^0(\widetilde{\mathfrak X},{\mathcal O}_{\widehat{\mathfrak X}}(k))$ and a certain space of ${\mathcal O}_{\widehat{K}}$-module valued harmonic cochains on the Bruhat-Tits tree of $G$. For $k\ge 2$ and even such an isomorphism was established by analytic methods in \cite{jer} whereas we proceed very algebro-geometrically in that we consequently reduce everything modulo the maximal ideal of ${\mathcal O}_{\widehat{K}}$ and work locally on the special fibre  $\widetilde{\mathfrak X}$. Finally, if $\kara(K)=0$, we demonstrate that integral structures are a strong tool for studying the "reduced"{}{} de Rham complex$${\mathcal R}_X^{\bullet}=[{\mathcal O}_X(-k)\stackrel{(\frac{d}{dz})^{k+1}}{\longrightarrow}{\mathcal O}_X(k+2)]$$on $X$ considered in \cite{schn}, \cite{ss}, for $k\ge0$ (here $z$ is a global variable on $X\subset\mathbb{P}^1_K$). It computes the de Rham cohomology $H^*(X,\Omega^{\bullet}_X\otimes\sym_K^k({\rm St}))$ of $X$ with coefficients in the $k$-th symmetric power $\sym_K^k({\rm St})$ of the standard representation of ${\rm SL}\sb 2(K)$. Its differential respects our integral structures, hence a complex$${\mathcal R}_{\widehat{\mathfrak X}}^{\bullet}=[{\mathcal O}_{\widehat{\mathfrak X}}(-k)\stackrel{(\frac{d}{dz})^{k+1}}{\longrightarrow}{\mathcal O}_{\widehat{\mathfrak X}}(k+2)]$$on ${\widehat{\mathfrak X}}$. We show that for $k>0$ we have $H^j(\widetilde{{\mathfrak X}},{\mathcal R}_{\widehat{{\mathfrak X}}}^{\bullet})=0$ for $j\ne 1$, while $H^1(\widetilde{{\mathfrak X}},{\mathcal R}_{\widehat{{\mathfrak X}}}^{\bullet})$ decomposes as $$H^1(\widetilde{{\mathfrak X}},{\mathcal R}_{\widehat{{\mathfrak X}}}^{\bullet})\cong H^1({\widetilde{\mathfrak{X}}},{\mathcal O}_{\widehat{\mathfrak{X}}}(-k))\oplus H^0(\widetilde{{\mathfrak X}},{\mathcal O}_{\widehat{{\mathfrak X}}}(k+2))\quad\quad\quad\quad(*).$$As an application, we show that structural features of the cohomology of varieties uniformized by $X$ can be deduced from $(*)$, thus show up already on $X$ (or rather $\widehat{\mathfrak X}$) itself. Namely we get the well known Hodge decomposition (first obtained by de Shalit \cite{deshhod}, see also \cite{schn}) $$H_{dR}^1(\Gamma\backslash X,\sym_K^k({\rm St}))=H^1(\Gamma,\sym_K^k({\rm St}))\oplus H^0(X_{\Gamma},{\mathcal O}_X(k+2)^{\Gamma})$$of $H_{dR}^1(X_{\Gamma},\sym_K^k({\rm St}))=H^1(\Gamma\backslash X,(\Omega_{{X}}^{\bullet}\otimes_{{K}}\sym_{{K}}^k({\rm St}))^{\Gamma})=H^1(\Gamma\backslash X,({\mathcal R}_X^{\bullet})^{\Gamma})$ simply by taking $\Gamma$-invariants for a cocompact discrete (torsionfree) subgroup $\Gamma<{\rm SL}\sb 2(K)$; no higher $\Gamma$-group cohomology is needed. Again, while earlier proofs were truly analytic we reduce everything to algebraic geometry on the irreducible components of $\widetilde{\mathfrak X}$ (these are all isomorphic to $\mathbb{P}_{\mathbb{F}}^1$). As a bonus of our method we obtain the degeneration of the "reduced"{} Hodge spectral sequence computing $H_{dR}^1(\Gamma\backslash X,\sym_K^k({\rm St}))$, as conjectured by Schneider \cite{schn}, and a complete description (in particular their dimensions) of the cohomology spaces $H^j(\Gamma\backslash X,{\mathcal O}_X(r)^{\Gamma})$ (any $j,r$). Moreover, for $k>0$, we describe a monodromy operator on $H^1(\widetilde{{\mathfrak X}},{\mathcal R}_{\widehat{{\mathfrak X}}}^{\bullet})$ as an isomorphism $H^0(\widetilde{{\mathfrak X}},{\mathcal O}_{\widehat{{\mathfrak X}}}(k+2))\cong H^1({\widetilde{\mathfrak{X}}},{\mathcal O}_{\widehat{\mathfrak{X}}}(-k))$. It induces the monodromy operator on $H_{dR}^1(\Gamma\backslash X,\sym_K^k({\rm St}))$ predicted by $p$-adic Hodge theory, so in particular we see that the latter respects integral de Rham structures (which in $p$-adic Hodge theory can not be expected in general) and that its monodromy filtration splits the Hodge filtration.

We mention that the integral structures in ${\mathcal O}_X(k)$ and in the "reduced"{} de Rham complex considered in this paper play an important role in the recent work of Breuil \cite{breuil}.

{\it Notations:} $K$ denotes a non-archimedean locally compact field and $K_a$ its algebraic closure, ${\cal O}_K$ its ring of integers, $\pi\in{\cal O}_K$ a fixed prime element and ${\mathbb{F}} $ the residue field with $q$ elements, $q\in p^{\mathbb{N}}$. We choose $\widehat{\pi}\in K_a$ such that $\widehat{\pi}^2=\pi$. Then $\widehat{K}=K(\widehat{\pi})$ is a ramified extension of $K$ of degree 2 with ring of integers ${\mathcal O}_{\widehat{K}}$. We let $\omega:K^{\times}_a\to\mathbb{Q}$ be the extension of the discrete valuation $\omega:K^{\times}\to\mathbb{Z}$ normalized by $\omega(\pi)=1$. For formal ${\mathcal O}_K$-schemes resp. $K$-rigid spaces we denote by a superscript $\widehat{.}$ the formal ${\mathcal O}_{\widehat{K}}$-schemes resp. $\widehat{K}$-rigid space obtained by the base change ${\mathcal O}_K\to{\mathcal O}_{\widehat{K}}$ resp. $K\to\widehat{K}$. 
For $E=K$ or $E=\widehat{K}$ and a formal (admissible) ${\cal O}_E$-scheme ${\mathfrak W}$ we let ${\mathfrak W}_E$ be its generic fibre, as a $E$-rigid space. We need the characters $\chi:G\to \widehat{K}^{\times}$, $\chi(\gamma)=\widehat{\pi}^{\omega(\det\gamma)}$, and $\varepsilon:G\to {\mathcal O}_K^{\times}$, $\varepsilon(\gamma)=\pi^{-\omega(\det\gamma)}\det\gamma$, of $G={\rm GL}\sb 2(K)$ and denote the Bruhat-Tits tree of $G$ by ${\mathcal {BT}}$. For $r\in\mathbb{R}$ we define $\lfloor r\rfloor, \lceil r\rceil\in\mathbb{Z}$ by requiring $\lfloor r\rfloor\le r<\lfloor r\rfloor+1$ and $\lceil r\rceil-1<r\le\lceil r\rceil$. 

\section{Integral structures in automorphic line bundles}
\label{fullgsec}

Let $X=\Omega^{(2)}_K$ be Drinfel'd's symmetric space of dimension $1$ over $K$. This is the $K$-rigid space obtained by removing all $K$-rational points from the projective line $\mathbb{P}^1_K$ over $K$. We choose a coordinate $z$ and define an action of $G$ on $X$ (on the left) by$$\gamma z=\frac{-b+az}{d-cz}\quad\mbox{for}\quad\gamma=\left(\begin{array}{*{2}c}a&b\\c&d\end{array}\right)\in G$$(\cite{jer} takes the other left action). Fix $k\in \mathbb{Z}$. For $f\in{\mathcal O}_{\widehat{X}}$ set\begin{gather}f|_{\gamma}(z)=\chi^k(\gamma)(a+cz)^{-k}f(\frac{b+dz}{a+cz})\quad\mbox{for}\quad\gamma=\left(\begin{array}{*{2}c}a&b\\c&d\end{array}\right)\in G.\label{xiii}\end{gather}
Denote by ${\mathcal O}_{\widehat{X}}(k)$ the structure sheaf of the $\widehat{K}$-rigid space $\widehat{X}$ endowed with the $G$-action on the left defined by (\ref{xiii}). (This is a {\it left} action; in \cite{jer} a right action is considered.)\\
As explained in \cite{jer}, the $K$-rigid space $X$ is the generic fibre of a certain $\pi$-adic strictly semistable formal ${\cal O}_K$-scheme $\mathfrak{X}$: the set $F^0$ of irreducible components of the reduction $\widetilde{\mathfrak{X}}$ of $\mathfrak{X}$ is in natural bijection with the set of vertices of ${\mathcal {BT}}$. Let $F^1$ be the set of subsets $\{Z_1,Z_2\}\subset F^0$ with $Z_1\cap Z_2\ne\emptyset$ and $Z_1\ne Z_2$; it corresponds to the set of edges of ${\mathcal {BT}}$. Each $Z\in F^0$ is isomorphic to $\mathbb{P}^1_{\mathbb{F}} $. The action of $G$ on $X$ extends to $\mathfrak{X}$. The admissible open subset$$U=\{P\in \mathbb{P}^1;\quad\omega(z(P))>-1\quad\mbox{and}\quad\omega(z(P)-x)<1\quad\mbox{for all}\quad x\in{\cal O}_K\}$$of $X$ is the tube (=preimage under the specialization map $X\to \mathfrak{X}$) of the central (with respect to $z$) irreducible component $Z_{\gamma_0}$ of $\widetilde{\mathfrak{X}}$. For $\gamma\in G$ define the irreducible component $Z_{\gamma}$ of $\widetilde{{\mathfrak{X}}}$ as $Z_{\gamma}=\gamma.Z_{\gamma_0}$. For $n\in\mathbb{Z}$ let $$\gamma_n=\left(\begin{array}{*{2}c}1&0\\0&\pi^n\end{array}\right)\in G.$$For a subset $E\subset F^0$ let $\widetilde{\mathfrak{U}}_E$ be the maximal open subscheme of $\widetilde{\mathfrak{X}}$ contained in $\cup_{Z\in E}Z$; in other words, the complement in $\widetilde{\mathfrak{X}}$ of the union of all irreducible components not in $E$. Let ${\mathfrak{U}}_E$ be the open formal subscheme of $\mathfrak{X}$ lifting $\widetilde{\mathfrak{U}}_E$. Letting $${\mathfrak{Y}}={\mathfrak{U}}_{\{Z_{\gamma_n};\,n\in\mathbb{Z}\}}$$ we have the open covering $${\mathfrak{X}}=\bigcup_{g\in {\rm SL}\sb 2(K)}g.{\mathfrak{Y}}$$ (${\rm SL}\sb 2(K)$ acts transitively on $F^1$). Let $f_{n,n}\in{\mathcal O}_{\mathfrak{Y}}(\mathfrak{U}_{\{Z_{\gamma_n},Z_{\gamma_{n+1}}\}})$ (resp. $f_{n,n+1}\in{\mathcal O}_{\mathfrak{Y}}(\mathfrak{U}_{\{Z_{\gamma_n},Z_{\gamma_{n+1}}\}})$) be an equation for the closed subscheme $Z_{\gamma_n}\cap \widetilde{\mathfrak{U}}_{\{Z_{\gamma_n},Z_{\gamma_{n+1}}\}}$ (resp. $Z_{\gamma_{n+1}}\cap \widetilde{\mathfrak{U}}_{\{Z_{\gamma_n},Z_{\gamma_{n+1}}\}}$) of $\mathfrak{U}_{\{Z_{\gamma_n},Z_{\gamma_{n+1}}\}}$. In local coordinates, there is an open embedding $\mathfrak{U}_{\{Z_{\gamma_n},Z_{\gamma_{n+1}}\}}\to\spf({\mathcal O}_K<X_1,X_2>/(X_1X_2-\pi))$ such that $f_{n,n}=X_1$ and $f_{n,n+1}=X_2$. Viewing $f_{n,n}$ and $f_{n,n+1}$ as sections of ${\mathcal O}_{\widehat{\mathfrak{Y}}}(\widehat{\mathfrak{U}}_{\{Z_{\gamma_n},Z_{\gamma_{n+1}}\}})$ we define$${\mathcal O}_{\widehat{\mathfrak U}_{\{Z_{\gamma_n},Z_{\gamma_{n+1}}\}}}(k)={\mathcal O}_{\widehat{\mathfrak U}_{\{Z_{\gamma_n},Z_{\gamma_{n+1}}\}}}.f_{n,n}^{\lceil\frac{kn}{2}\rceil}f_{n,n+1}^{\lceil\frac{k(n+1)}{2}\rceil}+{\mathcal O}_{\widehat{\mathfrak U}_{\{Z_{\gamma_n},Z_{\gamma_{n+1}}\}}}.\widehat{\pi}f_{n,n}^{\lfloor\frac{kn}{2}\rfloor}f_{n,n+1}^{\lfloor\frac{k(n+1)}{2}\rfloor},$$i.e. the ${\mathcal O}_{\widehat{\mathfrak{U}}_{\{Z_{\gamma_n},Z_{\gamma_{n+1}}\}}}$-submodule of ${\mathcal O}_{\widehat{\mathfrak{U}}_{\{Z_{\gamma_n},Z_{\gamma_{n+1}}\}}}\otimes_{{\mathcal O}_{\widehat{K}}}\widehat{K}$ generated by the two elements $f_{n,n}^{\lceil\frac{kn}{2}\rceil}f_{n,n+1}^{\lceil\frac{k(n+1)}{2}\rceil}$ and $\widehat{\pi}f_{n,n}^{\lfloor\frac{kn}{2}\rfloor}f_{n,n+1}^{\lfloor\frac{k(n+1)}{2}\rfloor}$. If $k$ is even this is just the line bundle generated by the element $z^{\frac{-k}{2}}$. If $k$ is odd this is not a line bundle; an explicit pair of generators is $\widehat{\pi}^{n+1}z^{\frac{-(k-1)}{2}}$, $\widehat{\pi}^{-n}z^{\frac{-(k+1)}{2}}$.

The ${\mathcal O}_{\widehat{\mathfrak{U}}_{\{Z_{\gamma_n},Z_{\gamma_{n+1}}\}}}(k)$ glue into an ${\mathcal O}_{\widehat{\mathfrak{Y}}}$-submodule ${\mathcal O}_{\widehat{\mathfrak{Y}}}(k)$ of ${\mathcal O}_{\widehat{\mathfrak{Y}}}\otimes_{{\mathcal O}_{\widehat{K}}}\widehat{K}$. Note that\begin{gather}{\mathcal O}_{\widehat{\mathfrak Y}}(k)|_{\widehat{\mathfrak U}_{\{Z_{\gamma_n}\}}}=\widehat{\pi}^{kn}{\mathcal O}_{\widehat{{\mathfrak U}}_{\{Z_{\gamma_n}\}}}\quad\mbox{inside}\quad{\mathcal O}_{\widehat{{\mathfrak U}}_{\{Z_{\gamma_n}\}}}\otimes_{{\mathcal O}_{\widehat{K}}}\widehat{K}.\label{xiv}\end{gather}As we remarked, if $k$ is even, ${\mathcal O}_{\widehat{\mathfrak{Y}}}(k)$ is the line bundle generated by the element $z^{\frac{-k}{2}}\in H^0(\widehat{\mathfrak{Y}},{\mathcal O}_{\widehat{\mathfrak{Y}}}\otimes_{{\mathcal O}_{\widehat{K}}}\widehat{K})$. For any $k$ again we have a canonical identification of sheaves $sp_*{\mathcal O}_{\widehat{X}}(k)={\mathcal O}_{\widehat{\mathfrak X}}\otimes_{{\mathcal O}_{\widehat{K}}}\widehat{K}$ where $sp:\widehat{X}\to\widehat{\mathfrak{X}}$ is the specialization map; we write $sp_*{\mathcal O}_{\widehat{X}}(k)$ when we refer to the $G$-equivariant structure on ${\mathcal O}_{\widehat{\mathfrak X}}\otimes_{{\mathcal O}_{\widehat{K}}}\widehat{K}$ induced by that on ${\mathcal O}_{\widehat{X}}(k)$.

\begin{pro}\label{gequivint} Let $\widehat{\mathfrak W}, \widehat{\mathfrak W}'$ be open formal subschemes of $\widehat{\mathfrak{Y}}$, let $\gamma\in G$ such that $\gamma\widehat{{\mathfrak W}}=\widehat{\mathfrak W}'$. Then the isomorphism$$\gamma:sp_*{\mathcal O}_{\widehat{X}}(k)|_{\widehat{\mathfrak W}}\cong sp_*{\mathcal O}_{\widehat{X}}(k)|_{\widehat{\mathfrak W}'}$$induces an isomorphism of subsheaves$$\gamma:{\mathcal O}_{\widehat{\mathfrak{Y}}}(k)|_{\widehat{\mathfrak W}}\cong{\mathcal O}_{\widehat{\mathfrak{Y}}}(k)|_{\widehat{\mathfrak W}'}.$$
\end{pro}

{\sc Proof:} (a) First we assume $\widehat{\mathfrak W}\subset \widehat{\mathfrak U}_{\{Z_{\gamma_n}\}}$ for some $n$; then also $\widehat{\mathfrak W}'\subset \widehat{\mathfrak U}_{\{Z_{\gamma_{n'}}\}}$ for some $n'$ and (\ref{xiv}) applies to $\widehat{\mathfrak W}$ and $\widehat{\mathfrak W}'$. In that situation we must show\begin{gather}2\omega((a+cz(P))^{-k})+k\omega(ad-bc)=k(n'-n)\quad\mbox{for}\quad\gamma=\left(\begin{array}{*{2}c}a&b\\c&d\end{array}\right)\label{xv}\end{gather} for each point $P$ in the generic fibre $\widehat{\mathfrak W}'_{\widehat{K}}$ of $\widehat{\mathfrak W}'$. Note that $\gamma Z_{\gamma_n}=Z_{\gamma_{n'}}$ and thus $\gamma_{n'}^{-1}\gamma\gamma_n$ stabilizes $Z_{\gamma_0}$, hence is an element of $K^{\times}.{\rm GL}\sb 2({\mathcal O}_K)$; in other words, $\gamma=\gamma_{n'}\delta\gamma_{n}^{-1}$ for some $\delta\in K^{\times}.{\rm GL}\sb 2({\mathcal O}_K)$. Therefore it suffices to check (\ref{xv}) in the cases\\(i) $\gamma=\gamma_m$ and $n'=n+m$ for some $m\in\mathbb{Z}$;\\(ii) $b=c=0=n=n'$ and $a=d$;\\(iii) $n=n'=0$ and $\gamma\in{\rm GL}\sb 2({\mathcal O}_K).$\\In either case (\ref{xv}) is immediate; for the case (iii) note that $\omega(z(P)-\beta)=0$ for any $\beta\in{\mathcal O}_{K}^{\times}$.

(b) Now let $\widehat{{\mathfrak W}}$, $\widehat{{\mathfrak W}}'$ be arbitrary. By construction, both ${\mathcal L}_1=\gamma_*({\mathcal O}_{\widehat{\mathfrak{Y}}}(k)|_{\widehat{\mathfrak W}})$ and ${\mathcal L}_2={\mathcal O}_{\widehat{\mathfrak{Y}}}(k)|_{\widehat{\mathfrak W}'}$ are ${\mathcal O}_{\widehat{\mathfrak W}'}$-modules contained in ${\mathcal O}_{\widehat{\mathfrak W}'}\otimes_{{\mathcal O}_{\widehat{K}}}\widehat{K}$ as lattices, i.e. ${\mathcal L}_i\otimes_{{\mathcal O}_{\widehat{K}}}\widehat{K}={\mathcal O}_{\widehat{\mathfrak W}'}\otimes_{{\mathcal O}_{\widehat{K}}}\widehat{K}$. By (a) we have ${\mathcal L}_1|_{\widehat{\mathfrak V}}={\mathcal L}_2|_{\widehat{\mathfrak V}}$ for an open formal subscheme ${\widehat{\mathfrak V}}$ of $\widehat{{\mathfrak W}}'$ whose reduction is dense in the reduction of $\widehat{\mathfrak W}'$. All this implies ${\mathcal L}_1={\mathcal L}_2$, using the following fact: for open formal subschemes $\widehat{{\mathfrak V}}_1\subset\widehat{{\mathfrak V}}_2$ of $\widehat{\mathfrak{Y}}$ with $\widehat{{\mathfrak V}}_1$ dense in $\widehat{{\mathfrak V}}_2$, and for $f\in({\mathcal O}_{\widehat{\mathfrak{Y}}}(k)\otimes_{{\mathcal O}_{\widehat{K}}}\widehat{K})(\widehat{{\mathfrak V}}_2)$ we have $f\in{\mathcal O}_{\widehat{\mathfrak{Y}}}(k)(\widehat{{\mathfrak V}}_2)$ if and only $f\in{\mathcal O}_{\widehat{\mathfrak{Y}}}(k)(\widehat{{\mathfrak V}}_1)$. To see this fact it suffices to show that for $g\in({\mathcal O}_{\widehat{\mathfrak{Y}}}(k)/(\widehat{\pi}))({\widehat{\mathfrak V}}_2)$ we have $g=0$ if and only $g|_{\widehat{{\mathfrak V}}_1}=0$ in $({\mathcal O}_{\widehat{\mathfrak{Y}}}(k)/(\widehat{\pi}))(\widehat{{\mathfrak V}}_1)$. This is immediate from the local analysis in section \ref{cosec} below.\hfill$\Box$\\

Thanks to \ref{gequivint} we can now move around ${\mathcal O}_{\widehat{\mathfrak Y}}(k)$ by means of the $G$-action on $\widehat{\mathfrak X}$ and obtain a $G$-equivariant coherent ${\mathcal O}_{\widehat{\mathfrak X}}$-module lattice ${\mathcal O}_{\widehat{\mathfrak X}}(k)$ inside $sp_*{\mathcal O}_{\widehat{X}}(k)$.

For $k_1, k_2\in\mathbb{Z}$ we have a $G$-equivariant surjective map (not needed in the sequel)$${\mathcal O}_{\widehat{\mathfrak X}}(k_1)\otimes_{{\mathcal O}_{\widehat{\mathfrak X}}} {\mathcal O}_{\widehat{\mathfrak X}}(k_2)\longrightarrow {\mathcal O}_{\widehat{\mathfrak X}}(k_1+k_2)$$which is multiplication of functions. This follows from equation (\ref{xiv}) and the argument in part (b) of the proof of \ref{gequivint}. It is as isomorphism if at least one of $k_1$ or $k_2$ is even, for in that case we are tensoring with a line bundle. On the other hand, it cannot be an isomorphism if both $k_1$ and $k_2$ are odd, because then the fibres of both ${\mathcal O}_{\widehat{\mathfrak X}}(k_j)$ at singular points of $\widetilde{X}$ are 2-dimensional, whereas ${\mathcal O}_{\widehat{\mathfrak X}}(k_1+k_2)$ is a line bundle (in this case).

\section{Cohomology}
\label{cosec}

For divisors $D$ on $\mathbb{P}_{\mathbb{F}} ^1$ let ${\mathcal L}(D)$ be the corresponding line bundle on ${\mathbb P}^1_{\mathbb{F}} $. By the usual convention, ${\mathcal L}(-D)\subset {\mathcal O}_{{\mathbb P}^1_{\mathbb{F}} }$ if $D$ is an effective divisor.
Fix a system $R$ of representatives for ${\mathbb{F}} $ in ${\mathcal O}_K$. For $a\in R$ and $n\in\mathbb{Z}$  let$$\gamma_{a,n}=\left(\begin{array}{*{2}c}1&\pi^{-n}a\\0&1\end{array}\right).$$An easy consideration on ${\mathcal {BT}}$ shows that$$\{Z_{\gamma_{n+1}}\}\bigcup\{Z_{\gamma_{a,n}\gamma_{n-1}};\,a\in R\}$$is the set of the $q+1$ many irreducible components of $\widetilde{\mathfrak X}$ meeting $Z_{\gamma_n}$. (The function $\pi^{n-1}z+\pi^{-1}a$ is a coordinate on ${\widehat{\mathfrak U}_{\{Z_{\gamma_{a,n}\gamma_{n-1}}\}}}$ in the sense that $\omega(\pi^{n-1}z(P)+\pi^{-1}a)=0$ for any $P\in({\widehat{\mathfrak U}_{\{Z_{\gamma_{a,n}\gamma_{n-1}}\}}})_{\widehat{K}}$.) Since $\gamma_{a,n}$ acts on $sp_*{\mathcal O}_{\widehat{X}}(k)$ with trivial automorphy factor it induces an isomorphism$$\gamma_{a,n}:\widehat{\pi}^{k(n-1)}{\mathcal O}_{\widehat{\mathfrak{X}}}(k)|_{\widehat{\mathfrak U}_{\{Z_{\gamma_{n-1}},Z_{\gamma_{n}}\}}}\cong\widehat{\pi}^{k(n-1)}{\mathcal O}_{\widehat{\mathfrak{X}}}(k)|_{\widehat{\mathfrak U}_{\{Z_{\gamma_{a,n}\gamma_{n-1}},Z_{\gamma_{n}}\}}}.$$Using this we can now give a local description of the $G$-equivariant coherent ${\mathcal O}_{\widetilde{\mathfrak X}}$-module ${\mathcal O}_{\widehat{\mathfrak X}}(k)/(\widehat{\pi})$ which we denote by ${\mathcal O}_{\widetilde{{\mathfrak X}}}(k)$.\\

{\it (a) First assume that $k$ is even}. Let $h_a\in{\mathcal O}_{Z_{\gamma_n}}$ be a local equation for $Z_{\gamma_n}\cap Z_{\gamma_{a,n}\gamma_{n-1}}$ in $Z_{\gamma_n}$, let $h_{\infty}\in{\mathcal O}_{Z_{\gamma_n}}$ be a local equation for $Z_{\gamma_n}\cap Z_{\gamma_{n+1}}$ in $Z_{\gamma_n}$. Then ${\mathcal O}_{\widetilde{\mathfrak X}}(k)\otimes_{{\mathcal O}_{\widetilde{\mathfrak X}}}{\mathcal O}_{Z_{\gamma_n}}$ is isomorphic to the following ${\mathcal O}_{Z_{\gamma_n}}$-submodule of the constant "rational function field"{} sheaf on $Z_{\gamma_n}\cong\mathbb{P}^1_{\mathbb{F}} $: locally around $Z_{\gamma_n}\cap Z_{\gamma_{a,n}\gamma_{n-1}}$ it is generated by $h_a^{\frac{k(n-1)}{2}-\frac{kn}{2}}$, locally around $Z_{\gamma_n}\cap Z_{\gamma_{n+1}}$ it is generated by $h_{\infty}^{\frac{k(n+1)}{2}-\frac{kn}{2}}$, and locally around other points it coincides with ${\mathcal O}_{Z_{\gamma_n}}$. Thus\begin{gather}{\mathcal O}_{\widetilde{\mathfrak X}}(k)\otimes_{{\mathcal O}_{\widetilde{\mathfrak X}}}{\mathcal O}_{Z}\cong{\mathcal L}(\frac{-k}{2}.\infty+\sum_{b\in {\mathbb{F}} }\frac{k}{2}.b)\label{vvv}\end{gather}for $Z=Z_{\gamma_n}$. By equivariance we get (\ref{vvv}) for any $Z\in F^0$. In particular, ${\mathcal O}_{\widetilde{\mathfrak X}}(k)\otimes_{{\mathcal O}_{\widetilde{\mathfrak X}}}{\mathcal O}_{Z}$ is of degree $\frac{(q-1)k}{2}$.

{\it (b) Now assume that $k$ is odd}. For $Z\in F^0$ let $$({\mathcal O}_{\widetilde{\mathfrak X}}(k)\otimes_{{\mathcal O}_{\widetilde{\mathfrak X}}}{\mathcal O}_Z)^c=\frac{{\mathcal O}_{\widetilde{\mathfrak X}}(k)\otimes_{{\mathcal O}_{\widetilde{\mathfrak X}}}{\mathcal O}_Z}{({\mathcal O}_{\widetilde{\mathfrak X}}(k)\otimes_{{\mathcal O}_{\widetilde{\mathfrak X}}}{\mathcal O}_Z)_{torsion}}.$$We then have$${\mathcal O}_{\widetilde{{\mathfrak X}}}(k)|_{\widetilde{\mathfrak{U}}_{\{Z_{\gamma_n},Z_{\gamma_{n+1}}\}}}=({\mathcal O}_{\widetilde{\mathfrak X}}(k)\otimes_{{\mathcal O}_{\widetilde{\mathfrak X}}}{\mathcal O}_{Z_{\gamma_{n}}})^c|_{\widetilde{\mathfrak{U}}_{\{Z_{\gamma_n},Z_{\gamma_{n+1}}\}}}\oplus ({\mathcal O}_{\widetilde{\mathfrak X}}(k)\otimes_{{\mathcal O}_{\widetilde{\mathfrak X}}}{\mathcal O}_{Z_{\gamma_{n+1}}})^c|_{\widetilde{\mathfrak{U}}_{\{Z_{\gamma_n},Z_{\gamma_{n+1}}\}}}.$$Explicitly, $({\mathcal O}_{\widetilde{\mathfrak X}}(k)\otimes_{{\mathcal O}_{\widetilde{\mathfrak X}}}{\mathcal O}_{Z_{\gamma_{n}}})^c|_{\widetilde{\mathfrak{U}}_{\{Z_{\gamma_n},Z_{\gamma_{n+1}}\}}}$ is generated by $f_{n,n}^{\frac{kn}{2}}f_{n,n+1}^{\frac{k(n+1)+1}{2}}$ if $n$ is even, and by $\widehat{\pi}f_{n,n}^{\frac{kn-1}{2}}f_{n,n+1}^{\frac{k(n+1)}{2}}$ if $n$ is odd. $({\mathcal O}_{\widetilde{\mathfrak X}}(k)\otimes_{{\mathcal O}_{\widetilde{\mathfrak X}}}{\mathcal O}_{Z_{\gamma_{n+1}}})^c|_{\widetilde{\mathfrak{U}}_{\{Z_{\gamma_n},Z_{\gamma_{n+1}}\}}}$ is generated by $\widehat{\pi}f_{n,n}^{\frac{kn}{2}}f_{n,n+1}^{\frac{k(n+1)-1}{2}}$ if $n$ is even, and by $f_{n,n}^{\frac{kn+1}{2}}f_{n,n+1}^{\frac{k(n+1)}{2}}$ if $n$ is odd. Now we proceed as in (a). By what we just saw, $({\mathcal O}_{\widetilde{\mathfrak X}}(k)\otimes_{{\mathcal O}_{\widetilde{\mathfrak X}}}{\mathcal O}_{Z_{\gamma_n}})^c$ is generated around $Z_{\gamma_{n+1}}\cap Z_{\gamma_n}$ by $h_{\infty}^{\lceil\frac{kn}{2}\rceil-\lceil\frac{k(n+1)}{2}\rceil}$, and around $Z_{\gamma_n}\cap Z_{\gamma_{a,n}\gamma_{n-1}}$ by $h_a^{\lceil\frac{k(n-1)}{2}\rceil-\lceil\frac{kn}{2}\rceil}$ (by equivariance, it suffices to check the latter for $a=0$). Thus\begin{gather}{\mathcal O}_{\widetilde{\mathfrak X}}(k)=\prod_{Z\in F^0}({\mathcal O}_{\widetilde{\mathfrak X}}(k)\otimes_{{\cal O}_{\widetilde{\mathfrak{X}}}}{\mathcal O}_Z)^c,\label{xvii}\end{gather} \begin{gather}({\mathcal O}_{\widetilde{\mathfrak X}}(k)\otimes_{{\mathcal O}_{\widetilde{\mathfrak X}}}{\mathcal O}_Z)^c\cong{\mathcal L}(\frac{-k-1}{2}.\infty+\sum_{b\in {\mathbb{F}} }\frac{k-1}{2}.b)\label{xvi}\end{gather}for $Z={Z_{\gamma_n}}$. By equivariance we get (\ref{xvi}) for any $Z\in F^0$. In particular, $({\mathcal O}_{\widetilde{\mathfrak X}}(k)\otimes_{{\mathcal O}_{\widetilde{\mathfrak X}}}{\mathcal O}_{Z})^c$ is of degree $\frac{(q-1)(k-1)}{2}-1$.

\begin{satz}\label{coho} (a) $H^*(\widetilde{\mathfrak{X}},{\mathcal O}_{\widehat{\mathfrak{X}}}(k))$ is ${\mathcal O}_{\widehat{K}}$-flat and $$H^*({\widetilde{\mathfrak{X}}},{\mathcal O}_{\widetilde{\mathfrak{X}}}(k))=H^*(\widetilde{\mathfrak{X}},{\mathcal O}_{\widehat{\mathfrak{X}}}(k))/(\pi).$$(b) For $k\le -1$ and also for $k=1$ we have $H^0(\widetilde{{\mathfrak{X}}},{\mathcal O}_{\widehat{\mathfrak{X}}}(k))=0$.\\(c) For $k\ge 0$ we have $H^1(\widetilde{{\mathfrak{X}}},{\mathcal O}_{\widehat{\mathfrak{X}}}(k))=0$. \end{satz}

{\sc Proof:} (i) First assume $k$ is even. To prove (c) it is enough to prove \begin{gather}\mathbb{R}^1\lim_{\stackrel{\leftarrow}{t}}H^0(\widetilde{{\mathfrak{X}}},{\mathcal O}_{\widehat{\mathfrak{X}}}(k)/(\widehat{\pi}^t))=0\label{v}\\\lim_{\stackrel{\leftarrow}{t}}H^1(\widetilde{{\mathfrak{X}}},{\mathcal O}_{\widehat{\mathfrak{X}}}(k)/(\widehat{\pi}^t))=0\label{vi}.\end{gather}For (\ref{v}) it suffices to show surjectivity of all transition maps $H^0(\widetilde{{\mathfrak{X}}},{\mathcal O}_{\widehat{\mathfrak{X}}}(k)/(\widehat{\pi}^{t+1}))\to H^0(\widetilde{{\mathfrak{X}}},{\mathcal O}_{\widehat{\mathfrak{X}}}(k)/(\widehat{\pi}^t))$. Using the long exact cohomology sequence associated with\begin{gather}0\to{\mathcal O}_{\widetilde{\mathfrak{X}}}(k)\stackrel{\widehat{\pi}^t}{\longrightarrow}{\mathcal O}_{\widehat{\mathfrak{X}}}(k)/(\widehat{\pi}^{t+1})\longrightarrow{\mathcal O}_{\widehat{\mathfrak{X}}}(k)/(\widehat{\pi}^t)\longrightarrow0\label{vii}\end{gather}this will be implied by\begin{gather}H^1(\widetilde{\mathfrak{X}},{\mathcal O}_{\widetilde{\mathfrak{X}}}(k))=0\label{viii}.\end{gather}Also (\ref{vi}) is reduced to (\ref{viii}) using (\ref{vii}), so let us prove (\ref{viii}). We have an exact sequence$$0\longrightarrow{\mathcal O}_{\widetilde{\mathfrak{X}}}(k)\longrightarrow\prod_{Z\in F^0}{\mathcal O}_{\widetilde{\mathfrak{X}}}(k)\otimes_{{\mathcal O}_{\widetilde{\mathfrak{X}}}} {\mathcal O}_{Z}\longrightarrow\prod_{\{Z_1,Z_2\}\in F^1}{\mathcal O}_{\widetilde{\mathfrak{X}}}(k)\otimes_{{\mathcal O}_{\widetilde{\mathfrak{X}}}} {\mathcal O}_{Z_1\cap Z_2}\longrightarrow0$$and a corresponding long exact sequence in cohomology. We know $$H^1(\widetilde{\mathfrak{X}},\prod_{Z\in F^0}{\mathcal O}_{\widetilde{\mathfrak{X}}}(k)\otimes_{{\mathcal O}_{\widetilde{\mathfrak{X}}}} {\mathcal O}_{Z})=\prod_{Z\in F^0}H^1(\widetilde{\mathfrak{X}},{\mathcal O}_{\widetilde{\mathfrak{X}}}(k)\otimes_{{\mathcal O}_{\widetilde{\mathfrak{X}}}} {\mathcal O}_{Z})=0$$because ${\mathcal O}_{\widetilde{\mathfrak{X}}}(k)\otimes_{{\mathcal O}_{\widetilde{\mathfrak{X}}}} {\mathcal O}_{Z}$ is isomorphic to a line bundle on $\mathbb{P}^1_{\mathbb{F}} \cong Z$ of non-negative degree as we saw above (since $k\ge 0$). On the other hand$$H^0(\widetilde{\mathfrak{X}},\prod_{Z\in F^0}{\mathcal O}_{\widetilde{\mathfrak{X}}}(k)\otimes_{{\mathcal O}_{\widetilde{\mathfrak{X}}}} {\mathcal O}_{Z})\longrightarrow H^0(\widetilde{\mathfrak{X}},\prod_{\{Z_1,Z_2\}\in F^1}{\mathcal O}_{\widetilde{\mathfrak{X}}}(k)\otimes_{{\mathcal O}_{\widetilde{\mathfrak{X}}}} {\mathcal O}_{Z_1\cap Z_2})$$is surjective: This follows from the contractiblity of ${\mathcal {BT}}$ and again the fact that each ${\mathcal O}_{\widetilde{\mathfrak{X}}}(k)\otimes_{{\mathcal O}_{\widetilde{\mathfrak{X}}}} {\mathcal O}_{Z}$ has non-negative degree, which implies that$$H^0(\widetilde{\mathfrak{X}},{\mathcal O}_{\widetilde{\mathfrak{X}}}(k)\otimes_{{\mathcal O}_{\widetilde{\mathfrak{X}}}} {\mathcal O}_{Z_1})\to H^0(\widetilde{\mathfrak{X}},{\mathcal O}_{\widetilde{\mathfrak{X}}}(k)\otimes_{{\mathcal O}_{\widetilde{\mathfrak{X}}}} {\mathcal O}_{Z_1\cap Z_2})$$for any $\{Z_1, Z_2\}\in F^1$ is surjective. To prove (b), since $H^0({\widetilde{\mathfrak{X}}},{\mathcal O}_{\widehat{\mathfrak{X}}}(k))=\lim_{\stackrel{\leftarrow}{t}}H^0(\widetilde{{\mathfrak{X}}},{\mathcal O}_{\widehat{\mathfrak{X}}}(k)/(\pi^t))$ we can reduce, using the long exact cohomolgy sequence associated with (\ref{vii}), to the statement$$H^0(\widetilde{{\mathfrak{X}}},{\mathcal O}_{\widetilde{\mathfrak{X}}}(k))=0.$$But this follows immediately from the injectivity of$$H^0(\widetilde{\mathfrak{X}},{\mathcal O}_{\widetilde{\mathfrak{X}}}(k))\longrightarrow H^0(\widetilde{\mathfrak{X}},\prod_{Z\in F^0}{\mathcal O}_{\widetilde{\mathfrak{X}}}(k)\otimes_{{\mathcal O}_{\widetilde{\mathfrak{X}}}} {\mathcal O}_{Z})$$and the fact that ${\mathcal O}_{\widetilde{\mathfrak{X}}}(k)\otimes_{{\mathcal O}_{\widetilde{\mathfrak{X}}}}{\mathcal O}_{Z}$ for each $Z\in F^0$ is isomorphic to a line bundle on $\mathbb{P}^1_{\mathbb{F}} \cong Z$ of negative degree as we saw above. To see the ${\mathcal O}_{\widehat{K}}$-flatness of $H^*(\widetilde{{\mathfrak{X}}},{\mathcal O}_{\widehat{\mathfrak{X}}}(k))$ in (a) we need to show injectivity of multiplication with $\widehat{\pi}$. This follows from (the proof of) (b) and (c) and the long exact cohomology sequence associated with$$0\longrightarrow {\mathcal O}_{\widehat{\mathfrak{X}}}(k)\stackrel{\widehat{\pi}}{\longrightarrow}{\mathcal O}_{\widehat{\mathfrak{X}}}(k)\longrightarrow{\mathcal O}_{\widetilde{\mathfrak{X}}}(k)\longrightarrow0.$$The base change statement follows similarly.\\
(ii) For odd $k$ the proofs are similar but easier in view of the decomposition (\ref{xvii}).\hfill$\Box$\\

The important vanishing $H^1(\widetilde{{\mathfrak{X}}},{\mathcal O}_{\widehat{\mathfrak{X}}}(k))=0$ was asserted for even $k\ge0$ in \cite{jer} Cor.24. However, the comparison with $H^1(\widehat{{{X}}},{\mathcal O}_{\widehat{{X}}}(k))$ invoked there does not seem to be justified.\\

Let $\Gamma<{\rm SL}\sb 2(K)$ be a cocompact discrete subgroup which for simplicity we assume to be torsion free (in general it contains a torsion free subgroup of finite index). Let $X_{\Gamma}=\Gamma\backslash X$, $\widehat{X}_{\Gamma}=\Gamma\backslash\widehat{X}$, $\widehat{\mathfrak X}_{\Gamma}=\Gamma\backslash\widehat{\mathfrak X}$ and $\widetilde{{\mathfrak X}}_{\Gamma}=\Gamma\backslash\widetilde{{\mathfrak X}}$ be the quotients for the free action by $\Gamma$; they all algebraize to projective schemes.

\begin{kor}\label{gammakoh} (a) For $k>0$ we have $$H^0(\widetilde{\mathfrak X}_{\Gamma},{\mathcal O}_{\widehat{\mathfrak{X}}}(-k)^{\Gamma})=0=H^1(\widetilde{\mathfrak X}_{\Gamma},{\mathcal O}_{\widehat{\mathfrak{X}}}(k+2)^{\Gamma}).$$In particular, $H^0({X}_{\Gamma},{\mathcal O}_{{{X}}}(-k)^{\Gamma})=0=H^1({X}_{\Gamma},{\mathcal O}_{{{X}}}(k+2)^{\Gamma})$.\\(b) $H^0(\widetilde{\mathfrak X}_{\Gamma},{\mathcal O}_{\widehat{\mathfrak{X}}}(k+2)^{\Gamma})$ and $H^1(\widetilde{\mathfrak X}_{\Gamma},{\mathcal O}_{\widehat{\mathfrak{X}}}(-k)^{\Gamma})$ are ${\mathcal O}_{\widehat{K}}$-flat and $$H^0(\widetilde{\mathfrak X}_{\Gamma},{\mathcal O}_{\widehat{\mathfrak{X}}}(k+2)^{\Gamma})\otimes_{{\mathcal O}_{\widehat{K}}}{\mathbb{F}} =H^0({\widetilde{\mathfrak X}}_{\Gamma},{\mathcal O}_{\widetilde{\mathfrak{X}}}(k+2)^{\Gamma})$$$$H^1(\widetilde{\mathfrak X}_{\Gamma},{\mathcal O}_{\widehat{\mathfrak{X}}}(-k)^{\Gamma})\otimes_{{\mathcal O}_{\widehat{K}}}{\mathbb{F}} =H^1({\widetilde{\mathfrak X}}_{\Gamma},{\mathcal O}_{\widetilde{\mathfrak{X}}}(-k)^{\Gamma}).$$(c) Serre duality identifies $H^1({X}_{\Gamma},{\mathcal O}_{{{X}}}(-k)^{\Gamma})$ with the dual of $H^0({X}_{\Gamma},{\mathcal O}_{{{X}}}(k+2)^{\Gamma})$.\\(d)$$H^j(\widetilde{\mathfrak{X}}_{\Gamma},{\mathcal O}_{\widehat{\mathfrak{X}}}(1)^{\Gamma})=0$$for any $j$. In particular, $H^j({{X}}_{\Gamma},{\mathcal O}_{{X}}(1)^{\Gamma})=0.$
\end{kor} 

{\sc Proof:} (a) For odd $k$ literally the same proof as in \ref{coho} applies, because in that case we have the decomposition (\ref{xvii}) which allows us to reduce to problems on each irreducible component --- these are the same for $\widetilde{\mathfrak{X}}$ and $\widetilde{\mathfrak{X}}_{\Gamma}$. Now let $k$ be even. From \ref{coho} we get $H^0(\widetilde{\mathfrak{X}},{\mathcal O}_{\widehat{\mathfrak{X}}}(-k))=0$ and $H^0(\widetilde{\mathfrak{X}},{\mathcal O}_{\widetilde{\mathfrak{X}}}(-k))=0$. In particular $H^0(\widetilde{\mathfrak{X}}_{\Gamma},{\mathcal O}_{\widehat{\mathfrak{X}}}(-k)^{\Gamma})=H^0(\widetilde{{\mathfrak{X}}},{\mathcal O}_{\widehat{{\mathfrak{X}}}}(-k))^{\Gamma}=0$ and $H^0(\widetilde{\mathfrak{X}}_{\Gamma},{\mathcal O}_{\widetilde{\mathfrak{X}}}(-k)^{\Gamma})=H^0({\widetilde{\mathfrak X}},{\mathcal O}_{\widetilde{\mathfrak{X}}}(-k))^{\Gamma}=0$. Now we have a ${\rm SL}\sb 2(K)$-equivariant isomorphism ${\mathcal O}_{{\widehat{\mathfrak{X}}}}(2)\cong \Omega^1_{{\widehat{\mathfrak{X}}}}$ on ${\widehat{\mathfrak{X}}}$, where $\Omega^1_{{\widehat{\mathfrak{X}}}}$ is the sheaf of relative logarithmic differentials for the log smooth formal $\spf({\mathcal O}_{\widehat{K}})$-scheme ${\widehat{\mathfrak{X}}}$ (with respect to the pull back log structures from the canonical log structures on ${\mathfrak{X}}$ and $\spf({\mathcal O}_{{K}})$). Thus ${\mathcal O}_{{\widehat{\mathfrak{X}}}}(2)^{\Gamma}$ can be identified with the sheaf of relative logarithmic differentials for the log smooth projective $\spec({\mathcal O}_{\widehat{K}})$-scheme ${\widehat{\mathfrak{X}}}_{\Gamma}$. This is a dualizing sheaf by \cite{dera} ch.I, sect.2, where it is called the sheaf of ${\it regular}$ differentials (the generalization to general projective log schemes is \cite{tsuji} Theorem 2.21). Since ${\mathcal O}_{\widehat{\mathfrak{X}}}(k+2)^{\Gamma}=({\mathcal O}_{\widehat{\mathfrak{X}}}(-k)^{\Gamma})^{\otimes(-1)}\otimes{\mathcal O}_{\widehat{\mathfrak{X}}}(2)^{\Gamma}$ (note that since $k$ is even we are dealing with line bundles here) we get $H^1({\widetilde{\mathfrak{X}}}_{\Gamma},{\mathcal O}_{\widehat{\mathfrak{X}}}(k+2)^{\Gamma})=0$ by Serre duality. The same argument works for the sheaves ${\mathcal O}_{\widetilde{\mathfrak{X}}}(.)$. For (b) we may now proceed as in \ref{coho}. For (c) note that ${\mathcal O}_{X}(2)$ is ${\rm SL}\sb 2(K)$-equivariantly isomorphic with the sheaf $\Omega^1_X$ of differentials on $X$, hence ${\mathcal O}_{{{X}}}(k+2)^{\Gamma}\cong ({\mathcal O}_{{{X}}}(-k)^{\Gamma})^{\otimes(-1)}\otimes\Omega^1_{X^{\Gamma}}$ (for even $k$ we just saw the integral version in (a)). The statements in (d) follow immediately from \ref{coho}.\hfill$\Box$

The fact $H^0(X_{\Gamma},{\mathcal O}_X(1)^{\Gamma})=0$ ("there are no non zero automorphic forms for $\Gamma$ of weight one") was proven by analytic methods in \cite{ss} Cor.13. For the $K$-vector space dimensions of $H^1(X_{\Gamma},{\mathcal O}_X(-k)^{\Gamma})$ and of $H^0(X_{\Gamma},{\mathcal O}_X(k+2)^{\Gamma})$ see \ref{hodge} below.  

\section{Modular representations}

Denote by ${\mathcal I}\subset{\mathcal O}_{\widetilde{\mathfrak X}}$ the ideal sheaf of functions vanishing at the singular points of $\widetilde{\mathfrak X}$. For $k\in\mathbb{Z}$ and $i\ge0$ let$${\mathcal O}_{\widetilde{\mathfrak X}}(k)(i)={\mathcal O}_{\widetilde{\mathfrak X}}(k)\otimes_{{\mathcal O}_{\widetilde{\mathfrak X}}}{\mathcal I}^i.$$ Let $Z\in F^0$. If $k$ is odd we let $$({\mathcal O}_{\widetilde{\mathfrak{X}}}(k)(i)\otimes_{{\cal O}_{\widetilde{\mathfrak{X}}}}{{\cal O}_Z})^c=({\mathcal O}_{\widetilde{\mathfrak{X}}}(k)\otimes_{{\cal O}_{\widetilde{\mathfrak{X}}}}{{\cal O}_Z})^c\otimes_{{\mathcal O}_{\widetilde{\mathfrak X}}}{\mathcal I}^i.$$
To unify notations, if $k$ is even we let $({\mathcal O}_{\widetilde{\mathfrak{X}}}(k)\otimes_{{\cal O}_{\widetilde{\mathfrak{X}}}}{{\cal O}_Z})^c={\mathcal O}_{\widetilde{\mathfrak{X}}}(k)\otimes_{{\cal O}_{\widetilde{\mathfrak{X}}}}{{\cal O}_Z}$ and$$({\mathcal O}_{\widetilde{\mathfrak{X}}}(k)(i)\otimes_{{\cal O}_{\widetilde{\mathfrak{X}}}}{{\cal O}_Z})^c={\mathcal O}_{\widetilde{\mathfrak{X}}}(k)\otimes_{{\cal O}_{\widetilde{\mathfrak{X}}}}{{\cal O}_Z}\otimes_{{\mathcal O}_{\widetilde{\mathfrak X}}}{\mathcal I}^i,$$i.e. for even $k$ the outer $(.)^c$ is redundant. We have\begin{gather}{\mathcal O}_{\widetilde{\mathfrak X}}(k)(i)=\prod_{Z\in F^0}({\mathcal O}_{\widetilde{\mathfrak X}}(k)(i)\otimes_{{\cal O}_{\widetilde{\mathfrak{X}}}}{\mathcal O}_Z)^c.\label{rrrr}\end{gather}if $k$ is odd and $i\ge0$ arbitrary, and also if $k$ is even and $i>0$. In particular, for such $(k,i)$ we have for any $Z\in F^0$ the natural injection$$\iota_Z:({\mathcal O}_{\widetilde{\mathfrak X}}(k)(i)\otimes_{{\cal O}_{\widetilde{\mathfrak{X}}}}{\mathcal O}_Z)^c\longrightarrow{\mathcal O}_{\widetilde{\mathfrak X}}(k)(i)$$and the canonical projection map $$\rho_{Z}:{\mathcal O}_{\widetilde{\mathfrak X}}(k)(i)\longrightarrow ({\mathcal O}_{\widetilde{\mathfrak X}}(k)(i)\otimes_{{\mathcal O}_{\widetilde{\mathfrak X}}}{\mathcal O}_{Z})^c.$$We denote maps induced by $\iota_{Z}$ resp. $\rho_{Z}$ in cohomology again by $\iota_{Z}$ resp. $\rho_{Z}$. 

\begin{lem}\label{indglob} Suppose $i\ge0$ if $k$ is odd, or $i>0$ if $k$ is even. Then we have a canonical $G$-equivariant isomorphism$$H^*(\widetilde{\mathfrak X},{\mathcal O}_{\widetilde{\mathfrak X}}(k)(i))\cong \ind_{K^{\times}{\rm GL}\sb 2({\mathcal O}_K)}^{G}H^*(\widetilde{\mathfrak X},({\mathcal O}_{\widetilde{\mathfrak X}}(k)(i)\otimes_{{\mathcal O}_{\widetilde{\mathfrak X}}}{\mathcal O}_{Z_{\gamma_0}})^c)$$
\end{lem}

{\sc Proof:} By definition, $\ind_{K^{\times}{\rm GL}\sb 2({\mathcal O}_K)}^{G}H^*(\widetilde{\mathfrak X},({\mathcal O}_{\widetilde{\mathfrak X}}(k)(i)\otimes_{{\mathcal O}_{\widetilde{\mathfrak X}}}{\mathcal O}_{Z_{\gamma_0}})^c)$ is the space of locally constant functions $u:G\to H^*(\widetilde{\mathfrak X},({\mathcal O}_{\widetilde{\mathfrak X}}(k)(i)\otimes_{{\mathcal O}_{\widetilde{\mathfrak X}}}{\mathcal O}_{Z_{\gamma_0}})^c)$ which satisfy $u(\eta\gamma)=\eta(u(\gamma))$ for $\eta\in K^{\times}{\rm GL}\sb 2({\mathcal O}_K)$, $\gamma\in G$. The action of $G$ is by $(\gamma.u)(\gamma')=u(\gamma'\gamma)$. Note that $K^{\times}{\rm GL}\sb 2({\mathcal O}_K)$ is the stabilizer of $Z_{\gamma_0}$ in $G$. Let $S\subset G$ be a subset such that $\gamma\mapsto Z_{\gamma}$ is a bijection between $S$ and $F^0$. The desired map is$$f\mapsto[u:G\to H^*(\widetilde{\mathfrak X},({\mathcal O}_{\widetilde{\mathfrak X}}(k)(i)\otimes_{{\mathcal O}_{\widetilde{\mathfrak X}}}{\mathcal O}_{Z_{\gamma_0}})^c),\,\gamma\mapsto\rho_{Z_{\gamma_0}}(\gamma.f)].$$Its inverse is$$[u:G\to H^*(\widetilde{\mathfrak X},({\mathcal O}_{\widetilde{\mathfrak X}}(k)(i)\otimes_{{\mathcal O}_{\widetilde{\mathfrak X}}}{\mathcal O}_{Z_{\gamma_0}})^c)]\mapsto\sum_{\gamma\in S}\gamma(\iota_{Z_{\gamma_0}}(u(\gamma^{-1}))).$$Note that $\gamma(\iota_{Z_{\gamma_0}}(u(\gamma^{-1})))$ is supported only on $Z_{\gamma}$.\hfill$\Box$

For a commutative ring $A$ and integers $n, s$ with $n\ge0$ let us denote by $\sym_A^n({\rm St})[s]$ the free $A$-module of homogeneous polynomials $F(X,Y)$ of degree $n$ in the variables $X, Y$ with coefficients in $A$, together with its ${\rm GL}\sb 2(A)$-action$$\gamma.F(X,Y)=(ad-bc)^sF(dX+bY,cX+aY)\quad\mbox{for}\quad\gamma=\left(\begin{array}{*{2}c}a&b\\c&d\end{array}\right)\in{\rm GL}\sb 2(A).$$Now consider for $k\in\mathbb{Z}$ the action of ${\rm GL}\sb 2({\mathbb{F}} )$ on ${\mathbb{F}} (z)$ given by\begin{gather}f|_{\gamma}(z)=(\frac{1}{a+cz})^kf(\frac{b+dz}{a+cz})\quad\mbox{for}\quad\gamma=\left(\begin{array}{*{2}c}a&b\\c&d\end{array}\right)\label{ix}.\end{gather}We view ${\mathbb{F}} (z)$ as the function field of $\mathbb{P}_{\mathbb{F}} ^1=\spec({\mathbb{F}} [z])\cup\{\infty\}$ and will consider line bundles on $\mathbb{P}_{\mathbb{F}} ^1$ stable for (\ref{ix}). Let $i\ge0$.

\begin{lem}\label{symgeom} (a) Suppose $k$ is even and $t=\frac{(q-1)k}{2}-i(q+1)\ge0$. Then, as ${\rm GL}\sb 2({\mathbb{F}} )$-representations, $$\sym_{\mathbb{F}} ^t({\rm St})[{i-\frac{k}{2}}]\cong H^0(\mathbb{P}_{\mathbb{F}} ^1,{\mathcal L}(\sum_{b\in {\mathbb{F}} }(\frac{k}{2}-i).b-(\frac{k}{2}+i).\infty)).$$(b) Suppose $k$ is odd and $t=\frac{(q-1)k-(q+1)}{2}-i(q+1)\ge0$. Then, as ${\rm GL}\sb 2({\mathbb{F}} )$-representations,$$\sym_{\mathbb{F}} ^t({\rm St})[{i-\frac{k-1}{2}}]\cong H^0(\mathbb{P}_{\mathbb{F}} ^1,{\mathcal L}(\sum_{b\in {\mathbb{F}} }(\frac{k-1}{2}-i).b-(\frac{k+1}{2}+i).\infty)).$$
\end{lem}

{\sc Proof:} In (a) the map sends $X^rY^{t-r}$ to $z^r(z-z^q)^{i-\frac{k}{2}}$ for $0\le r\le t$. In (b) it sends $X^rY^{t-r}$ to $z^r(z-z^q)^{i-\frac{k-1}{2}}$.\hfill$\Box$

We view $\sym_{\mathbb{F}} ^n({\rm St})[s]$ as a ${\rm GL}\sb 2({\mathcal O}_K)$-representation via the canonical map ${\rm GL}\sb 2({\mathcal O}_K)\to {\rm GL}\sb 2({\mathbb{F}} )$, and we then extend the action further to an action by $K^{\times}{\rm GL}\sb 2({\mathcal O}_K)$ by sending $\pi\in K^{\times}$ (i.e. the diagonal matrix with both entries equal to $\pi=\widehat{\pi}^2$) to the identity.

\begin{satz} (a) Suppose $k$ is even, $i>0$ and $t=\frac{(q-1)k}{2}-i(q+1)\ge0$. Then we have a canonical $G$-equivariant isomorphism$$H^0(\widetilde{\mathfrak X},{\mathcal O}_{\widetilde{\mathfrak X}}(k)(i))\cong \ind_{K^{\times}{\rm GL}\sb 2({\mathcal O}_K)}^{G}\sym_{\mathbb{F}} ^t({\rm St})[{i}-\frac{k}{2}].$$(b) Suppose $k$ is odd, $i\ge0$ and $t=\frac{(q-1)(k-1)}{2}-1-i(q+1)\ge0$. Then we have a canonical $G$-equivariant isomorphism$$H^0({\widetilde{\mathfrak X}},{\mathcal O}_{\widetilde{{\mathfrak X}}}(k)(i))\cong \ind_{K^{\times}{\rm GL}\sb 2({\mathcal O}_K)}^{G}\sym_{\mathbb{F}} ^t({\rm St})[i-\frac{k-1}{2}].$$\end{satz}

{\sc Proof:} We lift the ${\rm GL}\sb 2({\mathbb{F}} )$-action on $H^0(\mathbb{P}_{\mathbb{F}} ^1,{\mathcal L}(\sum_{b\in {\mathbb{F}} }(\frac{k}{2}-i).b-(\frac{k}{2}+i).\infty))$ if $k$ is even, resp. on $H^0(\mathbb{P}_{\mathbb{F}} ^1,{\mathcal L}(\sum_{b\in {\mathbb{F}} }(\frac{k-1}{2}-i).b-(\frac{k+1}{2}+i).\infty))$ if $k$ is odd, to an action by $K^{\times}{\rm GL}\sb 2({\mathcal O}_K)$ in the same way as explained for $\sym_{\mathbb{F}} ^n({\rm St})[s]$. Identifying the reduction of the global variable $z$ with our projective coordinate $z$ on $Z_{\gamma_0}\cong\mathbb{P}_{\mathbb{F}} ^1$ we use (\ref{xvi}) and (\ref{vvv}) to get $K^{\times}{\rm GL}\sb 2({\mathcal O}_K)$-equivariant isomorphisms$$H^0(\mathbb{P}_{\mathbb{F}} ^1,{\mathcal L}(\sum_{b\in {\mathbb{F}} }(\frac{k}{2}-i).b-(\frac{k}{2}+i).\infty))=H^0(\widetilde{\mathfrak X},({\mathcal O}_{\widetilde{\mathfrak X}}(k)(i)\otimes_{{\mathcal O}_{\widetilde{\mathfrak X}}}{\mathcal O}_{Z_{\gamma_0}})^c)$$if $k$ is even, resp.$$H^0(\mathbb{P}_{\mathbb{F}} ^1,{\mathcal L}(\sum_{b\in {\mathbb{F}} }(\frac{k-1}{2}-i).b-(\frac{k+1}{2}+i).\infty))=H^0(\widetilde{\mathfrak X},({\mathcal O}_{\widetilde{\mathfrak X}}(k)(i)\otimes_{{\mathcal O}_{\widetilde{\mathfrak X}}}{\mathcal O}_{Z_{\gamma_0}})^c)$$if $k$ is odd, thus we conclude by \ref{indglob} and \ref{symgeom}.\hfill$\Box$

We can now filter the representation $H^0(\widetilde{\mathfrak X},{\mathcal O}_{\widetilde{\mathfrak X}}(k))$ and determine its subquotients. For $k$ odd, $i\ge0$ and $t=\frac{(q-1)(k-1)}{2}-1-i(q+1)\ge q+1$ we have $H^1(\widetilde{\mathfrak X},({\mathcal O}_{\widetilde{\mathfrak X}}(k)(i+1)\otimes_{{\mathcal O}_{\widetilde{\mathfrak X}}}{\mathcal O}_{Z_{\gamma_0}})^c)=0$ (use (\ref{xvi})), hence$$\frac{\sym_{\mathbb{F}} ^{t}({\rm St})[i-\frac{k-1}{2}]}{\sym_{\mathbb{F}} ^{t-(q+1)}({\rm St})[i+1-\frac{k-1}{2}]}\cong\frac{H^0(\widetilde{\mathfrak X},({\mathcal O}_{\widetilde{\mathfrak X}}(k)(i)\otimes_{{\mathcal O}_{\widetilde{\mathfrak X}}}{\mathcal O}_{Z_{\gamma_0}})^c)}{H^0(\widetilde{\mathfrak X},({\mathcal O}_{\widetilde{\mathfrak X}}(k)(i+1)\otimes_{{\mathcal O}_{\widetilde{\mathfrak X}}}{\mathcal O}_{Z_{\gamma_0}})^c)}$$is a representation of  ${\rm GL}\sb 2(\mathbb{F})$ on the $(q+1)$-dimensional ${\mathbb{F}} $-vector space with basis the ${\mathbb{F}} $-rational points of ${\mathbb P}_{\mathbb{F}} ^1$. Explicitly, this is the quotient$$\frac{\sym_{\mathbb{F}} ^{t}({\rm St})[i-\frac{k-1}{2}]}{<X^jY^{t-j}-X^{q+j-1}Y^{t-q-j+1};\,\,1\le j\le t-q>_{\mathbb{F}}}.$$One might ask for its composition series. For example, if $q=2$, $k=9$, $i=0$, $t=3$, then the class of $X^3+Y^3+X^2Y$ (= the class of $X^3+Y^3+XY^2$) in this quotient spans a ${\rm GL}\sb 2(\mathbb{F})$-stable line. The results for even $k$ are similar, with $i>0$ and $t=\frac{(q-1)k}{2}-i(q+1)\ge q+1$, see also \cite{jer}. For the last $i$, the one for which $q\ge t\ge0$, we get $\sym_{\mathbb{F}} ^{t}({\rm St})[i-\frac{k-1}{2}]$ (if $k$ is odd), resp. $\sym_{\mathbb{F}} ^{t}({\rm St})[i-\frac{k}{2}]$ (if $k$ is even). To complete the picture it remains to observe that for even $k\ge4$ we have\begin{gather}\frac{H^0(\widetilde{\mathfrak X},{\mathcal O}_{\widetilde{\mathfrak X}}(k))}{H^0(\widetilde{\mathfrak X},{\mathcal O}_{\widetilde{\mathfrak X}}(k)(1))}\cong \ind_{N}^{G}{\bf 1}\label{x}\end{gather} where $N\subset G$ denotes the stabilizer of an (arbitrary) non-oriented edge $\{Z_1,Z_2\}\in F^1$ and ${\bf 1}$ its trivial representation: use $H^1(\widetilde{\mathfrak X},{\mathcal O}_{\widetilde{\mathfrak X}}(k)(1))=0$.\\

As an application, if $q$ is odd, Teitelbaum \cite{jer} constructs modular forms mod $\widehat{\pi}$ of weight $q+1$ (in fact, elements of $H^0(\widetilde{\mathfrak{X}},{\mathcal O}_{\widetilde{\mathfrak X}}(q+1)(\frac{q+1}{2}-1))$) for the entire group ${\rm SL}\sb 2(K)$. Here we will do the same if $q$ is even. The action of ${\rm SL}\sb 2(K)$ on the set $F^0$ has two orbits: the orbit $F^0_{even}$ of $Z_{\gamma_0}\in F^0$ and the orbit $F^0_{odd}$ of $Z_{\gamma_1}\in F^0$. Choose subsets $S_{even}$ and $S_{odd}$ of ${\rm SL}\sb 2(K)$ such that $\gamma\mapsto Z_{\gamma}$ defines bijections $S_{even}\cong F^0_{even}$ and $S_{odd}\cong F^0_{odd}$. Recall that we fixed a coordinate $z$ on $X$. For any $\gamma\in {\rm SL}\sb 2(K)$ we get another function $z\circ\gamma$ on $X$.

\begin{satz} The $H^0({\widetilde{\mathfrak X}},{\mathcal O}_{\widetilde{{\mathfrak X}}}(q+1))$-elements$${b}^+_{q+1}=\sum_{\gamma\in S_{even}}(\iota_{Z_{\gamma}}\circ\rho_{Z_{\gamma}})((z\circ\gamma^{-1}-(z\circ\gamma^{-1})^q)^{-1})$$$${b}^-_{q+1}=\sum_{\gamma\in S_{odd}}(\iota_{Z_{\gamma}}\circ\rho_{Z_{\gamma}})((z\circ\gamma^{-1}-(z\circ\gamma^{-1})^q)^{-1})$$are invariant for ${\rm SL}\sb 2(K)$, and interchanged by $\left(\begin{array}{*{2}c}0&1\\\pi&0\end{array}\right)$.
\end{satz}\hfill$\Box$

Now let us look at the modular representations $H^1(\widetilde{\mathfrak X},{\mathcal O}_{\widetilde{\mathfrak X}}(k))$ for $k<0$. If $k$ is even we get from$$0\longrightarrow{\mathcal O}_{\widetilde{\mathfrak{X}}}(k)\longrightarrow\prod_{Z\in F^0}{\mathcal O}_{\widetilde{\mathfrak{X}}}(k)\otimes_{{\mathcal O}_{\widetilde{\mathfrak{X}}}}{\mathcal O}_{Z}\longrightarrow\prod_{\{Z_1,Z_2\}\in F^1}{\mathcal O}_{\widetilde{\mathfrak{X}}}(k)\otimes_{{\mathcal O}_{\widetilde{\mathfrak{X}}}}{\mathcal O}_{Z_1\cap Z_2}\longrightarrow0$$the exact sequence$$0\to H^0(\widetilde{\mathfrak X},\prod_{\{Z_1,Z_2\}\in F^1}{\mathcal O}_{\widetilde{\mathfrak{X}}}(k)\otimes_{{\mathcal O}_{\widetilde{\mathfrak{X}}}}{\mathcal O}_{Z_1\cap Z_2})\longrightarrow$$$$ H^1(\widetilde{\mathfrak X},{\mathcal O}_{\widetilde{\mathfrak X}}(k))\longrightarrow H^1(\widetilde{\mathfrak X},\prod_{Z\in F^0}{\mathcal O}_{\widetilde{\mathfrak{X}}}(k)\otimes_{{\mathcal O}_{\widetilde{\mathfrak{X}}}}{\mathcal O}_{Z})\longrightarrow0.$$Here $H^0(\widetilde{\mathfrak X},\prod_{\{Z_1,Z_2\}\in F^1}{\mathcal O}_{\widetilde{\mathfrak{X}}}(k)\otimes_{{\mathcal O}_{\widetilde{\mathfrak{X}}}}{\mathcal O}_{Z_1\cap Z_2})$ is as in (\ref{x}). If $k$ is odd things are easier because then we have$$ H^1(\widetilde{\mathfrak X},{\mathcal O}_{\widetilde{\mathfrak X}}(k))\cong H^1(\widetilde{\mathfrak X},\prod_{Z\in F^0}({\mathcal O}_{\widetilde{\mathfrak{X}}}(k)\otimes_{{\mathcal O}_{\widetilde{\mathfrak{X}}}}{\mathcal O}_{Z})^c).$$Thus for any $k$, even or odd, we need to understand $H^1(\widetilde{\mathfrak X},\prod_{Z\in F^0}({\mathcal O}_{\widetilde{\mathfrak{X}}}(k)\otimes_{{\mathcal O}_{\widetilde{\mathfrak{X}}}}{\mathcal O}_{Z})^c)$ as a $G$-representation; by (the proof of) \ref{indglob} this means understanding $H^1(\widetilde{\mathfrak X},({\mathcal O}_{\widetilde{\mathfrak{X}}}(k)\otimes_{{\mathcal O}_{\widetilde{\mathfrak{X}}}}{\mathcal O}_{Z_{\gamma_0}})^c)$ as a ${\rm GL}\sb 2({\mathbb{F}} )$-representation. By an explicit computation on ${\mathbb P}^1_{\mathbb{F}}$, using the formulas (\ref{vvv}) and (\ref{xvi}), we see that Serre duality yields a ${\rm GL}\sb 2({\mathbb{F}} )$-equivariant isomorphism$$H^1(\widetilde{\mathfrak X},({\mathcal O}_{\widetilde{\mathfrak{X}}}(k)\otimes_{{\mathcal O}_{\widetilde{\mathfrak{X}}}}{\mathcal O}_{Z_{\gamma_0}})^c)\cong\Hom_{\mathbb{F}} (H^0(\widetilde{\mathfrak X},({\mathcal O}_{\widetilde{\mathfrak{X}}}(-k+2)(1)\otimes_{{\mathcal O}_{\widetilde{\mathfrak{X}}}}{\mathcal O}_{Z_{\gamma_0}})^c),{\mathbb{F}} )$$if $k$ is even, resp.$$H^1(\widetilde{\mathfrak X},({\mathcal O}_{\widetilde{\mathfrak{X}}}(k)\otimes_{{\mathcal O}_{\widetilde{\mathfrak{X}}}}{\mathcal O}_{Z_{\gamma_0}})^c)\cong\Hom_{\mathbb{F}} (H^0(\widetilde{\mathfrak X},({\mathcal O}_{\widetilde{\mathfrak{X}}}(-k+2)\otimes_{{\mathcal O}_{\widetilde{\mathfrak{X}}}}{\mathcal O}_{Z_{\gamma_0}})^c),{\mathbb{F}} )$$ if $k$ is odd. T(he duals of t)hese representations have been determined above. For example, for odd $k<0$, setting $t=\frac{(q-1)(-k-1)}{2}-1$ we get a canonical $G$-equivariant isomorphism$$H^1({\widetilde{\mathfrak X}},{\mathcal O}_{\widetilde{{\mathfrak X}}}(k))\cong \ind_{K^{\times}{\rm GL}\sb 2({\mathcal O}_K)}^{G}\ho_{\mathbb{F}} (\sym_{\mathbb{F}} ^t({\rm St})[\frac{k-1}{2}],{\mathbb{F}} ).$$
On the other hand, in section \ref{redhodge} below we will obtain for any $k<0$, even or odd, $G$-equivariant isomorphisms$$H^1({\widetilde{\mathfrak{X}}},{\mathcal O}_{\widetilde{\mathfrak{X}}}(k))\otimes\varepsilon^{-k-1}\cong H^0(\widetilde{{\mathfrak X}},{\mathcal O}_{\widetilde{{\mathfrak X}}}(2-k)).$$

\section{Harmonic cochains}
\label{harsec}

Fix $k\ge0$. On $\Hom_{\widehat{K}}(\sym_{\widehat{K}}^{k}({\rm St})[1]\otimes\chi^{-k-2},\widehat{K})$ the $G$-action is given by $(\gamma.h)(x)=h(\gamma^{-1}.x)$ for $\gamma\in G$, $x\in \sym_{\widehat{K}}^{k}({\rm St})[1]\otimes\chi^{-k-2}$ and $h\in\Hom_{\widehat{K}}(\sym_{\widehat{K}}^{k}({\rm St})[1]\otimes\chi^{-k-2},{\widehat{K}})$. --- (In everything here and below we could replace $\Hom_{\widehat{K}}(\sym_{\widehat{K}}^{k}({\rm St})[1]\otimes\chi^{-k-2},\widehat{K})$ by the isomorphic $G$-representation $\sym_{\widehat{K}}^{k}({\rm St})[-k-1]\otimes\chi^{k+2}$: the isomorphism sends $h_j$ (as defined below) to $X^{k-j}Y^{j}$.) --- We set$$C^1(k+2)=\prod_{\{Z_1,Z_2\}\in{F^1}}\Hom_{\widehat{K}}(\sym_{\widehat{K}}^{k}({\rm St})[1]\otimes\chi^{-k-2},\widehat{K}),$$
$$C^0(k+2)=\prod_{Z\in F^0}\Hom_{\widehat{K}}(\sym_{\widehat{K}}^{k}({\rm St})[1]\otimes\chi^{-k-2},\widehat{K})$$(products of copies of $\Hom_{\widehat{K}}(\sym_{\widehat{K}}^{k}({\rm St})[1]\otimes\chi^{-k-2},\widehat{K})$, indexed by $F^1$ resp. $F^0$). On $C^1(k+2)$ we define a $G$-action by$$(\gamma.f)_{\{Z_1,Z_2\}}=\gamma(f_{\gamma^{-1}\{Z_1,Z_2\}})$$for $\gamma\in G$ and $(f_{\{Z_1,Z_2\}})_{\{Z_1,Z_2\}}\in C^1(k+2)$. For $Z\in F^0$ let ${\rm sg}(Z)=1$ if $Z\in F^0_{even}$ and ${\rm sg}(Z)=-1$ if $Z\in F^0_{odd}$. Moreover let $$*(Z)=\{Z'\in F^0;\,\{Z,Z'\}\in F^1\}.$$Then we have the operator$$C^1(k+2)\stackrel{\Delta}{\longrightarrow}C^0(k+2),\quad\quad(f_{\{Z_1,Z_2\}})_{\{Z_1,Z_2\}}\mapsto({\rm sg}(Z)\sum_{Z'\in*(Z)}f_{\{Z,Z'\}})_{Z}$$and we define $C^1_{har}(k+2)$ by the exact sequence$$0\longrightarrow C^1_{har}(k+2)\longrightarrow C^1(k+2)\stackrel{\Delta}{\longrightarrow}C^0(k+2).$$

This is the variant with non-trivial coefficients of the space $C^1_{har}(\widehat{K})$ of $\widehat{K}$-valued harmonic cochains on ${\mathcal {BT}}$ which is defined by the exact sequence \begin{gather}0\longrightarrow C^1_{har}(\widehat{K})\longrightarrow \prod_{\{Z_1,Z_2\}\in{F^1}}\widehat{K}\stackrel{\Delta}{\longrightarrow}\prod_{Z\in F^0}\widehat{K}.\end{gather}Let $\Omega^1_{{\widehat{\mathfrak{X}}}}$ denote the sheaf of logarithmic differential forms for the morphism of log schemes $\widehat{{\mathfrak X}}\to\spf({\cal O}_{\widehat{K}})$ (with log structures defined by the respective special fibres). Define$$res:\Gamma(\widehat{{\mathfrak X}},\Omega^1_{{\widehat{\mathfrak{X}}}})\to C^1_{har}(\widehat{K})$$to be the unique $G$-equivariant morphism of ${\cal O}_{\widehat{K}}$-modules with$$res(\eta)_{\{Z_{\gamma_0},Z_{\gamma_{-1}}\}}=a_{-1}$$for $\eta\in\Gamma(\widehat{{\mathfrak X}},\Omega^1_{{\widehat{\mathfrak{X}}}})$, where $$\eta(z)=\sum_{j\in\mathbb{Z}}a_jz^jdz$$ is the Laurent expansion of $\eta$ on the annulus $]Z_{\gamma_0}\cap Z_{\gamma_{-1}}[=sp^{-1}(Z_{\gamma_0}\cap Z_{\gamma_{-1}})\subset \widehat{X}$ reducing to $Z_{\gamma_0}\cap Z_{\gamma_{-1}}$. (That $res(\eta)$ indeed lies in $C^1_{har}(\widehat{K})$ follows from the residue theorem on ${\mathbb P}^1$.) This map also has a version with non-trivial coefficients, as follows. Consider the $G$-equivariant map$$\Gamma(\widehat{{\mathfrak X}},{\mathcal O}_{\widehat{{\mathfrak X}}}(k+2))\longrightarrow\Hom_{\widehat{K}}(\sym_{\widehat{K}}^{k}({\rm St})[1]\otimes\chi^{-k-2}, \Gamma(\widehat{{\mathfrak X}},\Omega^1_{{\widehat{\mathfrak{X}}}})),\quad g\mapsto\Phi_g$$where $\Phi_g$ is defined by$$\Phi_g(X^iY^{k-i})=g(z)z^idz,\quad\quad0\le i\le k.$$We use it to define the $G$-equivariant map$$Res^0:\Gamma(\widehat{{\mathfrak X}},{\mathcal O}_{\widehat{{\mathfrak X}}}(k+2))\to \Hom(\sym_{\widehat{K}}^{k}({\rm St})[1]\otimes\chi^{-k-2},C^1_{har}(\widehat{K}))=C_{har}^1(k+2)$$$$g\mapsto res\circ\phi_g.$$We will work with the following more explicit description of $Res^0$: it is the unique $G$-equivariant morphism of ${\cal O}_{\widehat{K}}$-modules with$$(Res^0(g)_{\{Z_{\gamma_0},Z_{\gamma_{-1}}\}})(X^iY^{k-i})=a_{-i-1}$$for $g\in\Gamma(\widehat{{\mathfrak X}},{\mathcal O}_{\widehat{{\mathfrak X}}}(k+2))$ and $0\le i\le k$, where $$g(z)=\sum_{j\in\mathbb{Z}}a_jz^j$$ is the Laurent expansion of $g$ on the annulus $]Z_{\gamma_0}\cap Z_{\gamma_{-1}}[=sp^{-1}(Z_{\gamma_0}\cap Z_{\gamma_{-1}})\subset \widehat{X}$ reducing to $Z_{\gamma_0}\cap Z_{\gamma_{-1}}$. Equivalently, $(Res^0(g)_{\{Z_{1},Z_{2}\}})(X^iY^{k-i})$ for arbitrary $\{Z_1,Z_2\}\in{F^1}$ can be described follows. Choose a $\gamma\in G$ such that $\gamma.\{Z_{1},Z_{2}\}=\{Z_{\gamma_0},Z_{\gamma_{-1}}\}$. Let $\sum_{j\in\mathbb{Z}}a_jz^j$ be the Laurent expansion of $\gamma.g$ on $]Z_{\gamma_0}\cap Z_{\gamma_{-1}}[$ and write $$\gamma.(X^iY^{k-i})=\sum_{s=0}^kc_sX^sY^{k-s}$$ in $\sym_{\widehat{K}}^{k}({\rm St})[1]\otimes\chi^{-k-2}$. Then $(Res^0(g)_{\{Z_{1},Z_{2}\}})(X^iY^{k-i})=\sum_{s=0}^ka_{-s-1}c_s$. This is independent on the choice of $\gamma$.

We want to show that $Res^0$ is injective and to describe its image. For $Z\in F^0$ choose $\gamma\in G$ with $Z=Z_{\gamma}$ and define$$L_{Z}=\gamma.\Hom_{{\cal O}_{\widehat{K}}}(\sym_{{\cal O}_{\widehat{K}}}^{k}({\rm St})[1],{\cal O}_{\widehat{K}})\subset\Hom_{\widehat{K}}(\sym_{\widehat{K}}^{k}({\rm St})[1]\otimes\chi^{-k-2},\widehat{K}).$$In this definition we consider $\Hom_{{\cal O}_{\widehat{K}}}(\sym_{{\cal O}_{\widehat{K}}}^{k}({\rm St})[1],{\cal O}_{\widehat{K}})$ not as a ${\rm GL}\sb 2({\cal O}_{\widehat{K}})$-representation but only as a ${\cal O}_{\widehat{K}}$-submodule of the $\widehat{K}$-vector space underlying the $G$-representation $\Hom_{\widehat{K}}(\sym_{\widehat{K}}^{k}({\rm St})[1]\otimes\chi^{-k-2},\widehat{K})$. For $\{Z_1,Z_2\}\in{F^1}$ we write $L_{\{Z_1,Z_2\}}=L_{Z_1}\cap L_{Z_2}$ and then let$$Z^1(k+2)=\prod_{\{Z_1,Z_2\}\in{F^1}}L_{\{Z_1,Z_2\}},$$$$Z^0(k+2)=\prod_{Z}\in{F^0}L_{Z},$$subspaces of $C^1(k+2)$ resp. of $C^0(k+2)$. We define $Z^1_{har}(k+2)$ by the exact sequence\begin{gather}0\longrightarrow Z^1_{har}(k+2)\longrightarrow Z^1(k+2)\stackrel{\Delta}{\longrightarrow}\prod_{Z\in F^0}L_Z.\label{koketex}\end{gather}

\begin{lem} The image of $Res^0$ lies in $Z^1_{har}(k+2)$.
\end{lem}

{\sc Proof:} By $G$-equivariance it suffices to check $Res^0(g)_{\{Z_{\gamma_0},Z_{\gamma_{-1}}\}}\in L_{Z_{\gamma_0}}\cap L_{Z_{\gamma_{-1}}}$ for all $g\in\Gamma(\widehat{{\mathfrak X}},{\mathcal O}_{\widehat{{\mathfrak X}}}(k+2))$. Let $g(z)=\sum_{j\in\mathbb{Z}}a_jz^j$ be the Laurent expansion of $g$ on $]Z_{\gamma_0}\cap Z_{\gamma_{-1}}[$. From (\ref{xiv}) we deduce \begin{gather}\omega(g(P))\ge 0\quad\mbox{for all closed points }\, P\in]\widetilde{\mathfrak{U}}_{\{Z_{\gamma_0}\}}[\label{ineq1}\end{gather} \begin{gather}\omega(g(P))\ge\frac{-k-2}{2}\quad\mbox{for all closed points }\, P\in]\widetilde{\mathfrak{U}}_{\{Z_{\gamma_{-1}}\}}[.\label{ineq2}\end{gather} From (\ref{ineq1}) we get $\omega(a_j)\ge0$ for all $j$ (with a point $P\in ]Z_{\gamma_0}\cap Z_{\gamma_{-1}}[$ approach $]\widetilde{\mathfrak{U}}_{\{Z_{\gamma_0}\}}[$), hence $Res^0(g)_{\{Z_{\gamma_0},Z_{\gamma_{-1}}\}}\in \Hom_{{\cal O}_{\widehat{K}}}(\sym_{{\cal O}_{\widehat{K}}}^{k}({\rm St})[1],{\cal O}_{\widehat{K}})=L_{Z_{\gamma_0}}$. From (\ref{ineq2}) we get $\omega(a_j)\ge\frac{-(k+2)-2j}{2}$ for all $j$ (with a point $P\in ]Z_{\gamma_0}\cap Z_{\gamma_{-1}}[$ approach $]\widetilde{\mathfrak{U}}_{\{Z_{\gamma_{-1}}\}}[$). Now in $\sym_{\widehat{K}}^{k}({\rm St})[1]\otimes\chi^{-k-2}$ we have $\gamma_{-1}(X^iY^{k-i})=\widehat{\pi}^{k-2i}X^iY^{k-i}$. Thus $Res^0(g)_{\{Z_{\gamma_0},Z_{\gamma_{-1}}\}}(\gamma_{-1}(X^iY^{k-i}))=\widehat{\pi}^{k-2i}Res^0(g)_{\{Z_{\gamma_0},Z_{\gamma_{-1}}\}}(X^iY^{k-i})=\widehat{\pi}^{k-2i}a_{-i-1}$ lies in ${\cal O}_{\widehat{K}}$, thus $\gamma_1.Res^0(g)_{\{Z_{\gamma_0},Z_{\gamma_{-1}}\}}$ lies in $\Hom_{{\cal O}_{\widehat{K}}}(\sym_{{\cal O}_{\widehat{K}}}^{k}({\rm St})[1],{\cal O}_{\widehat{K}})$, thus $Res^0(g)_{\{Z_{\gamma_0},Z_{\gamma_{-1}}\}}$ lies in $L_{Z_{\gamma_{-1}}}$.

\begin{satz}\label{harm} $$Res^0:\Gamma(\widehat{{\mathfrak X}},{\mathcal O}_{\widehat{{\mathfrak X}}}(k+2))\longrightarrow Z^1_{har}(k+2)$$ is an isomorphism.
\end{satz}

{\sc Proof:} (i) First we claim that the sequence (\ref{koketex}) is also exact on the right. Let $\widetilde{Z}^1(k+2)=Z^1(k+2)/(\widehat{\pi})$ and for $Z\in F^0$ let $\widetilde{L}_{Z}=L_{Z}/(\widehat{\pi})$. Then it is enough to show that the map$$\widetilde{Z}^1(k+2)\stackrel{\widetilde{\Delta}}{\longrightarrow}\prod_{Z\in F^0}\widetilde{L}_Z$$induced by $\Delta$ is surjective. For $\{Z_1, Z_2\}\in{F^1}$ let $$D_{\{Z_1,Z_2\}}^{Z_1}=\bi(L_{\{Z_1,Z_2\}}\to\widetilde{L}_{Z_1})$$$$E_{\{Z_1,Z_2\}}=\bi(L_{\{Z_1,Z_2\}}\to({L}_{Z_1}+L_{Z_2})/(\widehat{\pi}))$$(images under the natural maps). Note that $\dim_{\mathbb{F}} (D_{\{Z_1,Z_2\}}^{Z_1})=\frac{k+2}{2}$ and $\dim_{\mathbb{F}} (E_{\{Z_1,Z_2\}})=1$ if $k$ is even, and $\dim_{\mathbb{F}} (D_{\{Z_1,Z_2\}}^{Z_1})=\frac{k+1}{2}$ and $E_{\{Z_1,Z_2\}}=0$ if $k$ is odd (for explicit descriptions see below). For $Z\in F^0$ let$$\widetilde{Z}^1(k+2)_Z=\prod_{Z'\in*(Z)}D_{\{Z,Z'\}}^Z.$$Then $\widetilde{\Delta}$ factors as$$\widetilde{Z}^1(k+2)\stackrel{\beta}{\longrightarrow}\prod_{Z\in F^0}\widetilde{Z}^1(k+2)_Z\stackrel{\delta=\prod\delta_Z}{\longrightarrow}\prod_{Z\in F^0}\widetilde{L}_Z$$where $\beta$ is the product of the natural projection maps. We have an exact sequence$$0\longrightarrow\widetilde{Z}^1(k+2)\stackrel{\beta}{\longrightarrow}\prod_{Z\in F^0}\widetilde{Z}^1(k+2)_Z\stackrel{\alpha}{\longrightarrow}\prod_{\{Z_1,Z_2\}}E_{\{Z_1,Z_2\}}$$where $\alpha$ is defined as$$\alpha(((g_{Z,Z'})_{Z'\in*(Z)})_{Z\in F^0})_{\{Z_1,Z_2\}}={\rm sg}(Z_1)g_{Z_1,Z_2}+{\rm sg}(Z_2)g_{Z_2,Z_1}$$for $\{Z_1,Z_2\}\in F^1$. For $Z\in F^0$ we define $\widetilde{Z}^1_{har}(k+2)_Z$ by the exact sequence$$0\longrightarrow\widetilde{Z}^1_{har}(k+2)_Z\stackrel{\nu_Z}{\longrightarrow}\widetilde{Z}^1(k+2)_Z\stackrel{\delta_Z}{\longrightarrow}\widetilde{L}_Z.$$Now it is enough to prove that each $\delta_Z$ (and hence $\delta$) is surjective, and that$$\prod_{Z\in F^0}\widetilde{Z}^1_{har}(k+2)_Z\stackrel{\alpha\circ(\prod_Z\nu_Z)}{\longrightarrow}\prod_{\{Z_1,Z_2\}}E_{\{Z_1,Z_2\}}$$is surjective. The surjectivity of $\alpha\circ(\prod_Z\nu_Z)$, an empty statement if $k$ is odd, will be implied by the surjectivity of its factors$$\widetilde{Z}^1_{har}(k+2)_{Z_1}\stackrel{\mu_{Z_1,Z_2}}{\longrightarrow}E_{\{Z_1,Z_2\}}.$$ Let us make the objects explicit. By equivariance we may assume $Z=Z_{\gamma_0}$, resp. $\{Z_1,Z_2\}=\{Z_{\gamma_0}, Z_{\gamma_{-1}}\}$. For $0\le j\le k$ define $h_j\in\Hom_{\widehat{K}}(\sym_{\widehat{K}}^{k}({\rm St})[1]\otimes\chi^{-k-2},{\widehat{K}})$ by
$$h_j(X^iY^{k-i})=\left\{\begin{array}{l@{\quad:\quad}l}1&\quad i=j\\0&\quad i\ne j\end{array}\right..$$Then one finds $$L_{\gamma_0}=\oplus_{j=0}^k{\mathcal O}_{\widehat{K}}.h_j,\quad L_{\gamma_{1}}=\oplus_{j=0}^k(\widehat{\pi}^{k-2j}).h_j,\quad\quad L_{\gamma_{-1}}=\oplus_{j=0}^k(\widehat{\pi}^{2j-k}).h_j,$$$$D_{\{Z_{\gamma_{-1}},Z_{\gamma_0}\}}^{Z_{\gamma_0}}=\oplus_{j=0}^{\lfloor\frac{k}{2}\rfloor}{\mathbb{F}} .h_j,\quad\quad D_{\{Z_{\gamma_{1}},Z_{\gamma_0}\}}^{Z_{\gamma_0}}=\oplus_{j=\lceil\frac{k}{2}\rceil}^{k}{\mathbb{F}} .h_j$$and if $k$ is even also $E_{\{Z_{\gamma_{-1}},Z_{\gamma_0}\}}={\mathbb{F}} .h_{\frac{k}{2}}$. The surjectivity of $\delta_{Z_{\gamma_0}}$ follows from $\widetilde{L}_{Z_{\gamma_0}}= D_{\{Z_{\gamma_{-1}},Z_{\gamma_0}\}}^{Z_{\gamma_0}}+D_{\{Z_{\gamma_{1}},Z_{\gamma_0}\}}^{Z_{\gamma_0}}$. For the surjectivity of $\mu_{Z_{\gamma_{0}},Z_{\gamma_{-1}}}$ (if $k$ is even): the element $h_{\frac{k}{2}}\in E_{\{Z_{\gamma_{-1}},Z_{\gamma_0}\}}$ is the image of the $\widetilde{Z}^1_{har}(k+2)_{Z_{\gamma_0}}$-element with entry $h_{\frac{k}{2}}$ in the $\{Z_{\gamma_{-1}},Z_{\gamma_0}\}$-component, with entry $-h_{\frac{k}{2}}$ in the $\{Z_{\gamma_{1}},Z_{\gamma_0}\}$-component, and with entry $0$ at all other components.\\
(ii) Let $\widetilde{Z}_{har}^1(k+2)=Z_{har}^1(k+2)/(\widehat{\pi})$. To prove the theorem, since $\Gamma(\widehat{{\mathfrak X}},{\mathcal O}_{\widehat{{\mathfrak X}}}(k+2))$ and $Z_{har}^1(k+2)$ are $\widehat{\pi}$-adically complete and separated, and since $Z_{har}^1(k+2)$ is ${\cal O}_{\widehat{K}}$-flat, it is enough to prove that the induced map$$\widetilde{Res}^0:\Gamma(\widehat{{\mathfrak X}},{\mathcal O}_{\widehat{{\mathfrak X}}}(k+2))/(\widehat{\pi})\to \widetilde{Z}_{har}^1(k+2)$$is an isomorphism. Since $\prod_{Z\in F^0}{L}_Z$ is ${\cal O}_{\widehat{K}}$-flat it follows from (i) that (\ref{koketex}) reduces modulo $(\widehat{\pi})$ to an exact sequence $$0\longrightarrow\widetilde{Z}_{har}^1(k+2)\longrightarrow\widetilde{Z}^1(k+2)\stackrel{\widetilde{\Delta}}{\longrightarrow}\prod_{Z\in F^0}\widetilde{L}_Z.$$We then also obtain from (i) for any $Z\in F^0$ exact sequences$$0\to \widetilde{Z}_{har}^1(k+2)\to \prod_{Z\in F^0}\widetilde{Z}^1_{har}(k+2)_Z\to\prod_{\{Z_1,Z_2\}\in F^1}E_{\{Z_1,Z_2\}}\to0$$and (by surjectivity of $\delta_Z$) the estimates$$\dim_{\mathbb{F}} (\widetilde{Z}^1_{har}(k+2)_Z)=\left\{\begin{array}{l@{\quad:\quad}l}\frac{(q-1)(k+1)}{2}&k\mbox{ odd}\\\frac{(q-1)(k+2)}{2}+1&k\mbox{ even}\end{array}\right.$$Now let us look at the source of $\widetilde{Res}^0$. By \ref{coho} we know that this is $H^0({{\mathfrak{X}}},{\mathcal O}_{\widetilde{\mathfrak{X}}}(k+2))$. Our discussion in section \ref{cosec} implies that the natural restriction maps induce an exact sequence (note $H^1({{\mathfrak{X}}},{\mathcal O}_{\widetilde{\mathfrak{X}}}(k+2))=0$)$$0\to H^0({{\mathfrak{X}}},{\mathcal O}_{\widetilde{\mathfrak{X}}}(k+2))\to\prod_{Z\in F^0}H^0({{\mathfrak{X}}},({\mathcal O}_{\widetilde{\mathfrak{X}}}(k+2)\otimes_{{\cal O}_{\widetilde{\mathfrak{X}}}}{{\cal O}_Z})^c)\to\prod_{\{Z_1,Z_2\}}J_{\{Z_1,Z_2\}}\to0$$where $\dim_{\mathbb{F}} (J_{\{Z_1,Z_2\}})=1$ if $k$ is even, and $J_{\{Z_1,Z_2\}}=0$ if $k$ is odd. Since $\widetilde{Res}^0$ induces isomorphisms $J_{\{Z_1,Z_2\}}\cong E_{\{Z_1,Z_2\}}$ it now suffices to see that the map $$H^0({{\mathfrak{X}}},({\mathcal O}_{\widetilde{\mathfrak{X}}}(k+2)\otimes_{{\cal O}_{\widetilde{\mathfrak{X}}}}{{\cal O}_Z})^c)\to\widetilde{Z}^1_{har}(k+2)_Z$$induced by $\widetilde{Res}^0$ is an isomorphism for any $Z\in F^0$, or, by equivariance, for $Z=Z_{\gamma_0}$. Recall that identifying $\mathbb{P}^1_{\mathbb{F}} \cong Z_{\gamma_0}$ as before we have$$({\mathcal O}_{\widetilde{\mathfrak{X}}}(k+2)\otimes_{{\cal O}_{\widetilde{\mathfrak{X}}}}{{\cal O}_{Z_{\gamma_0}}})^c\cong\left\{\begin{array}{l@{\quad:\quad}l}{\mathcal L}(\frac{-k-3}{2}.\infty+\sum_{b\in {\mathbb{F}} }\frac{k+1}{2}.b)&k\mbox{ odd}\\{\mathcal L}(\frac{-k-2}{2}.\infty+\sum_{b\in {\mathbb{F}} }\frac{k+2}{2}.b)&k\mbox{ even}\end{array}\right..$$For $k$ odd, if $g\in H^0({{\mathfrak{X}}},({\mathcal O}_{\widetilde{\mathfrak{X}}}(k+2)\otimes_{{\cal O}_{\widetilde{\mathfrak{X}}}}{{\cal O}_{Z_{\gamma_0}}})^c)=H^0(\mathbb{P}^1_{\mathbb{F}} ,{\mathcal L}(\frac{-k-3}{2}.\infty+\sum_{b\in {\mathbb{F}} }\frac{k+1}{2}.b))$ lies in the kernel of $\widetilde{Res}^0$ then it is an element even of $H^0(\mathbb{P}^1_{\mathbb{F}} ,{\mathcal L}(\frac{-k-3}{2}.\infty))$ and therefore it vanishes. Thus $\widetilde{Res}^0$ is injective. Similarly for even $k$. On the other hand by our above computation we find $\dim_{\mathbb{F}} (H^0({{\mathfrak{X}}},({\mathcal O}_{\widetilde{\mathfrak{X}}}(k+2)\otimes_{{\cal O}_{\widetilde{\mathfrak{X}}}}{{\cal O}_{Z_{\gamma_0}}})^c))=\dim_{\mathbb{F}} (\widetilde{Z}^1_{har}(k+2)_Z)$, thus $\widetilde{Res}^0$ is also surjective and the proof is complete.\hfill$\Box$

The $p$-adic Shimura isomorphism \cite{ss} p.98 is an immediate consequence of \ref{harm}.

\section{The reduced de Rham complex}
\label{redhodge}

In this section $\kara(K)=0$. Fix $k\ge0$ and for our fixed coordinate $z$ let $\partial=\frac{d}{dz}$. Let $(\Omega_X^{\bullet}\otimes_K\sym_K^k({\rm St}),\partial\otimes\id)$ be the de Rham complex on $X$ with coefficients in $\sym_K^k({\rm St})$. By \cite{ss} p.97 this complex is ${\rm SL}\sb 2(K)$-equivariantly quasi-isomorphic with the "reduced de Rham complex"$${\mathcal R}_X^{\bullet}=[{\mathcal O}_X(-k)\stackrel{\partial^{k+1}}{\longrightarrow}{\mathcal O}_X(k+2)]$$on $X$. (The genesis of this "theta operator"{} $\partial^{k+1}$ from $(\Omega_X^{\bullet}\otimes_K\sym_K^k({\rm St}),\partial\otimes\id)$ is completely parallel to that of the theta operator on classical modular forms, cf. \cite{cole}).

We change bases $K\to\widehat{K}$. Since $\omega(z(P))=-n$ for any $n$ and any point $P\in(\widehat{\mathfrak U}_{\{Z_{\gamma_n}\}})_{\widehat{K}}$ the operator $\partial$ on ${\mathcal O}_{\widehat{{\mathfrak U}}_{\{Z_{\gamma_n}\}}}\otimes_{{\cal O}_{\widehat{K}}}{\widehat{K}}=sp_*{\mathcal O}_{(\widehat{{\mathfrak U}}_{\{Z_{\gamma_n}\}})_{\widehat{K}}}$ restricts to a map $\partial:{\mathcal O}_{\widehat{{\mathfrak U}}_{\{Z_{\gamma_n}\}}}\to\pi^n.{\mathcal O}_{\widehat{{\mathfrak U}}_{\{Z_{\gamma_n}\}}}$. Iterating we get a map $\partial^{k+1}:\widehat{\pi}^{-{kn}}.{\mathcal O}_{\widehat{{\mathfrak U}}_{\{Z_{\gamma_n}\}}}\to\widehat{\pi}^{{(k+2)n}}.{\mathcal O}_{{\mathfrak U}_{\{Z_{\gamma_n}\}}}$, i.e. a map $\partial^{k+1}:{\mathcal O}_{\widehat{\mathfrak X}}(-k)|_{\widehat{{\mathfrak U}}_{\{Z_{\gamma_n}\}}}\to{\mathcal O}_{\widehat{\mathfrak X}}(k+2)|_{\widehat{{\mathfrak U}}_{\{Z_{\gamma_n}\}}}$. By equivariance we see that $\partial^{k+1}$ induces a map $\partial^{k+1}:{\mathcal O}_{\widehat{\mathfrak X}}(-k)|_{\widehat{{\mathfrak U}}_{\{Z\}}}\to{\mathcal O}_{\widehat{\mathfrak X}}(k+2)|_{\widehat{{\mathfrak U}}_{\{Z\}}}$ for {\it any} $Z\in F^0$. By an argument similar to that at the end of the proof of \ref{gequivint} it follows that $\partial^{k+1}$ respects these integral structures also above the singular points of $\widetilde{\mathfrak{X}}$, hence a complex$${\mathcal R}_{\widehat{{\mathfrak X}}}^{\bullet}=[{\mathcal O}_{\widehat{{\mathfrak X}}}(-k)\stackrel{\partial^{k+1}}{\longrightarrow}{\mathcal O}_{\widehat{{\mathfrak X}}}(k+2)].$$We denote by ${\cal H}^i({\mathcal R}_{\widehat{{\mathfrak X}}}^{\bullet})$ for $i=0$ and $i=1$ the cohomology sheaves.

\begin{satz}\label{symdec} For any $i,j$ we have canonical isomorphisms$$H^j(\widetilde{{\mathfrak X}},{\cal H}^i({\mathcal R}_{\widehat{{\mathfrak X}}}^{\bullet}))\cong H^j(\widetilde{{\mathfrak X}},{\mathcal R}_{\widehat{{\mathfrak X}}}^{i}).$$
\end{satz}

{\sc Proof:} For $i=0$ the map is induced by the canonical injection ${\cal H}^0({\mathcal R}_{\widehat{{\mathfrak X}}}^{\bullet})\to{\mathcal R}_{\widehat{{\mathfrak X}}}^{0}={\mathcal O}_{\widehat{{\mathfrak X}}}(-k)$, for $i=1$ it is induced by the canonical surjection ${\mathcal R}_{\widehat{{\mathfrak X}}}^{1}\to{\cal H}^1({\mathcal R}_{\widehat{{\mathfrak X}}}^{\bullet})$. Once we know the claim for $i=0$ it follows that $H^*(\widetilde{{\mathfrak X}},{\cal B})=0$ for ${\cal B}=\bi({\mathcal R}_{\widehat{{\mathfrak X}}}^{0}\to{\mathcal R}_{\widehat{{\mathfrak X}}}^{1})=\ke({\mathcal R}_{\widehat{{\mathfrak X}}}^{1}\to{\cal H}^1({\mathcal R}_{\widehat{{\mathfrak X}}}^{\bullet}))$, hence the claim for $i=1$. Thus we concentrate on the case $i=0$. Denote by $(.)_m$ reduction modulo $\widehat{\pi}^m$. Since ${\cal H}^0({\mathcal R}_{\widehat{{\mathfrak X}}}^{\bullet})=\lim_{\stackrel{\leftarrow}{m}}({\cal H}^0({\mathcal R}_{\widehat{{\mathfrak X}}}^{\bullet}))_m$ and ${\mathcal R}_{\widehat{{\mathfrak X}}}^{0}=\lim_{\stackrel{\leftarrow}{m}}({\mathcal R}_{\widehat{{\mathfrak X}}}^{0})_m$, the spectral sequence for the composition of derived functors $\mathbb{R}\lim_{\stackrel{\leftarrow}{m}}\mathbb{R}\Gamma(\widehat{{\mathfrak X}},.)$ shows that it suffices to show$$H^j(\widetilde{{\mathfrak X}}, ({\cal H}^0({\mathcal R}_{\widehat{{\mathfrak X}}}^{\bullet}))_m)\cong H^j(\widetilde{{\mathfrak X}},({\mathcal R}_{\widehat{{\mathfrak X}}}^{0})_m)$$for any $m$. Now ${\mathcal R}_{\widehat{{\mathfrak X}}}^{0}$ and hence also its subsheaf ${\cal H}^0({\mathcal R}_{\widehat{{\mathfrak X}}}^{\bullet})$ is ${\mathcal O}_{\widehat{K}}$-flat. Therefore one gets exact sequences of sheaves$$0\to{\mathcal F}_{m-1}\stackrel{\widehat{\pi}^{m-1}}{\longrightarrow}{\mathcal F}_m\to{\mathcal F}_1\to0$$ for ${\mathcal F}={\mathcal R}_{\widehat{{\mathfrak X}}}^{0}$ and ${\mathcal F}={\cal H}^0({\mathcal R}_{\widehat{{\mathfrak X}}}^{\bullet})$. Using the associated long exact cohomology sequences we reduce our task to proving the isomorphism just stated in the case $m=1$. Now observe that ${\cal H}^0({\mathcal R}_{{X}}^{\bullet}\otimes_K\widehat{K})$ is precisely the locally constant sheaf generated by the $\widehat{K}$-vector space of polynomials in the variable $z$ of degree at most $k$. Thus ${\cal H}^0({\mathcal R}_{\widehat{{\mathfrak X}}}^{\bullet})$ consists of such polynomials subject to growth conditions. Namely, since ${\mathcal R}_{\widehat{{\mathfrak X}}}^{0}|_{\widehat{\mathfrak{U}}_{\{Z_{\gamma_n}\}}}={\cal O}_{\widehat{\mathfrak{X}}}(-k)|_{\widehat{\mathfrak{U}}_{\{Z_{\gamma_n}\}}}=\widehat{\pi}^{-kn}{\cal O}_{\widehat{\mathfrak{X}}}|_{\widehat{\mathfrak{U}}_{\{Z_{\gamma_n}\}}}$ and $\omega(z(P))=-n$ for any $n$ and any point $P\in(\widehat{\mathfrak U}_{\{Z_{\gamma_n}\}})_{\widehat{K}}$ we have$${\cal H}^0({\mathcal R}_{\widehat{{\mathfrak X}}}^{\bullet})(\widehat{\mathfrak{U}}_{\{Z_{\gamma_n}\}})=\{\sum_{0\le t\le k}d_tz^t|\quad d_t\in\widehat{K},\,\omega(d_t)\ge tn-\frac{kn}{2}\},$$\begin{align}{\cal H}^0({\mathcal R}_{\widehat{{\mathfrak X}}}^{\bullet})(\widehat{\mathfrak{U}}_{\{Z_{\gamma_n},Z_{\gamma_{n-1}}\}})&={\cal H}^0({\mathcal R}_{\widehat{{\mathfrak X}}}^{\bullet})(\widehat{\mathfrak{U}}_{\{Z_{\gamma_n}\}})\cap{\cal H}^0({\mathcal R}_{\widehat{{\mathfrak X}}}^{\bullet})(\widehat{\mathfrak{U}}_{\{Z_{\gamma_{n-1}}\}})\notag\\{}&=\{\sum_{0\le t\le k}d_tz^t|\quad d_t\in\widehat{K},\,\omega(d_t)\ge\left\{\begin{array}{l@{\quad:\quad}l}tn-\frac{kn}{2}&t\ge\frac{k}{2}\\t(n-1)-\frac{k(n-1)}{2}&t\le\frac{k}{2}\end{array}\right.\}\notag\end{align}(any $k$, even or odd). For $Z\in F^0$ let $({\cal H}^0({\mathcal R}_{\widehat{{\mathfrak X}}}^{\bullet}))_1^Z$ be the image of the composition$${\cal H}^0({\mathcal R}_{\widehat{{\mathfrak X}}}^{\bullet})\to{\mathcal R}_{\widehat{{\mathfrak X}}}^{0}={\mathcal O}_{\widehat{{\mathfrak X}}}(-k)\to({\cal O}_{\widetilde{\mathfrak{X}}}(-k)\otimes_{{\mathcal O}_{\widetilde{{\mathfrak X}}}}{\mathcal O}_Z)^c.$$ Then the above shows$$({\cal H}^0({\mathcal R}_{\widehat{{\mathfrak X}}}^{\bullet}))_1^{Z_{\gamma_n}}(\widehat{\mathfrak{U}}_{\{Z_{\gamma_n}\}})=\{\sum_{0\le t\le k}\overline{d}_tz^t|\quad \overline{d}_t\in\frac{(\widehat{\pi}^{(2t-k)n})}{(\widehat{\pi}^{(2t-k)n+1})}\}$$$$({\cal H}^0({\mathcal R}_{\widehat{{\mathfrak X}}}^{\bullet}))_1^{Z_{\gamma_n}}(\widehat{\mathfrak{U}}_{\{Z_{\gamma_n},Z_{\gamma_{n-1}}\}})=\{\sum_{\frac{k}{2}\le t\le k}\overline{d}_tz^t|\quad \overline{d}_t\in\frac{(\widehat{\pi}^{(2t-k)n})}{(\widehat{\pi}^{(2t-k)n+1})}\}$$$$({\cal H}^0({\mathcal R}_{\widehat{{\mathfrak X}}}^{\bullet}))_1^{Z_{\gamma_{n}}}(\widehat{\mathfrak{U}}_{\{Z_{\gamma_{n}},Z_{\gamma_{n+1}}\}})=\{\sum_{0\le t\le\frac{k}{2}}\overline{d}_tz^t|\quad \overline{d}_t\in\frac{(\widehat{\pi}^{(2t-k)n})}{(\widehat{\pi}^{(2t-k)n+1})}\}.$$Similar descriptions hold at other $Z\in F^0$, resp. $\{Z_1,Z_2\}\in F^1$, by equivariance. We find$$({\cal H}^0({\mathcal R}_{\widehat{{\mathfrak X}}}^{\bullet}))_1=\prod_{Z\in F^0}({\cal H}^0({\mathcal R}_{\widehat{{\mathfrak X}}}^{\bullet}))_1^Z$$if $k$ is odd (because then there are no summands $\overline{d}_{\frac{k}{2}}z^{\frac{k}{2}}$ to consider). If $k$ is even we find an exact sequence$$0\to({\cal H}^0({\mathcal R}_{\widehat{{\mathfrak X}}}^{\bullet}))_1\to\prod_{Z\in F^0}({\cal H}^0({\mathcal R}_{\widehat{{\mathfrak X}}}^{\bullet}))_1^Z\to\prod_{\{Z_1,Z_2\}\in F^1}({\cal H}^0({\mathcal R}_{\widehat{{\mathfrak X}}}^{\bullet}))_1^{Z_1,Z_2}\to0$$where $({\cal H}^0({\mathcal R}_{\widehat{{\mathfrak X}}}^{\bullet}))_1^{Z_1,Z_2}$ for $\{Z_1,Z_2\}\in F^1$ is a sheaf with $({\cal H}^0({\mathcal R}_{\widehat{{\mathfrak X}}}^{\bullet}))_1^{Z_1,Z_2}(U)\cong {\mathbb{F}} $ if $U\cap Z_1\cap Z_2\ne\emptyset$, and $=0$ for other open $U\subset \widetilde{{\mathfrak X}}$. On the other hand we have$$({\mathcal R}_{\widehat{{\mathfrak X}}}^{0})_1={\cal O}_{\widetilde{\mathfrak{X}}}(-k)=\prod_{Z\in F^0}({\cal O}_{\widetilde{\mathfrak{X}}}(-k)\otimes_{{\mathcal O}_{\widetilde{{\mathfrak X}}}}{\mathcal O}_Z)^c$$if $k$ is odd, and an exact sequence$$0\to({\mathcal R}_{\widehat{{\mathfrak X}}}^{0})_1\to\prod_{Z\in F^0}({\cal O}_{\widetilde{\mathfrak{X}}}(-k)\otimes_{{\mathcal O}_{\widetilde{{\mathfrak X}}}}{\mathcal O}_Z)^c\to\prod_{\{Z_1,Z_2\}\in F^1}{\cal O}_{\widetilde{\mathfrak{X}}}(-k)\otimes_{{\mathcal O}_{\widetilde{{\mathfrak X}}}}{\mathcal O}_{Z_1\cap Z_2}\to0$$if $k$ is even. For $Z\in F^0$ let$$\alpha_Z:({\cal H}^0({\mathcal R}_{\widehat{{\mathfrak X}}}^{\bullet}))_1^Z\to({\cal O}_{\widetilde{\mathfrak{X}}}(-k)\otimes_{{\mathcal O}_{\widetilde{{\mathfrak X}}}}{\mathcal O}_Z)^c$$be the inclusion. If $k$ is even then for $\{Z_1,Z_2\}\in F^1$ the maps $\alpha_{Z_1}$ and $\alpha_{Z_2}$ commute with obvious isomorphisms$$\alpha_{Z_1,Z_2}:({\cal H}^0({\mathcal R}_{\widehat{{\mathfrak X}}}^{\bullet}))_1^{Z_1,Z_2}\to{\cal O}_{\widetilde{\mathfrak{X}}}(-k)\otimes_{{\mathcal O}_{\widetilde{{\mathfrak X}}}}{\mathcal O}_{Z_1\cap Z_2}.$$Since the  $\alpha_Z$ also commute with our map $({\cal H}^0({\mathcal R}_{\widehat{{\mathfrak X}}}^{\bullet}))_1\to({\mathcal R}_{\widehat{{\mathfrak X}}}^{0})_1$ in question, it remains to prove that the $\alpha_Z$ induce isomorphisms in cohomology. By equivariance it is enough to do this for $Z=Z_{\gamma_0}$. We identify $\spec({\mathbb{F}} [z])\cup\{\infty\}=\mathbb{P}^1_{\mathbb{F}} \cong Z_{\gamma_0}$ such that this $z$ on $\mathbb{P}^1_{\mathbb{F}} $ is induced by the global variable $z$ on $X$. In particular, $\infty\in\mathbb{P}^1_{\mathbb{F}} $ corresponds to $Z_{\gamma_0}\cap Z_{\gamma_1}$, and $0\in\mathbb{P}^1_{\mathbb{F}} $ corresponds to $Z_{\gamma_0}\cap Z_{\gamma_{-1}}$. Let $\iota:\mathbb{P}^1_{\mathbb{F}} \cong Z_{\gamma_0}\to\widetilde{{\mathfrak X}}$ be the closed immersion. Since we have$$H^*(\widetilde{{\mathfrak X}},{\mathcal F})=H^*(\mathbb{P}^1_{\mathbb{F}} ,\iota^{-1}{\mathcal F})$$for both ${\mathcal F}=({\cal H}^0({\mathcal R}_{\widehat{{\mathfrak X}}}^{\bullet}))_1^{Z_{\gamma_0}}$ and ${\mathcal F}=({\cal O}_{\widetilde{\mathfrak{X}}}(-k)\otimes_{{\mathcal O}_{\widetilde{{\mathfrak X}}}}{\mathcal O}_{Z_{\gamma_0}})^c$, we must show that$$H^*(\mathbb{P}^1_{\mathbb{F}} ,\iota^{-1}({\cal H}^0({\mathcal R}_{\widehat{{\mathfrak X}}}^{\bullet}))_1^{Z_{\gamma_0}})\to H^*(\mathbb{P}^1_{\mathbb{F}} ,\iota^{-1}({\cal O}_{\widetilde{\mathfrak{X}}}(-k)\otimes_{{\mathcal O}_{\widetilde{{\mathfrak X}}}}{\mathcal O}_{Z_{\gamma_0}})^c)$$is an isomorphism. If on $\mathbb{P}^1_{\mathbb{F}} $ we define the divisor$$D=\left\{\begin{array}{l@{\quad:\quad}l}\frac{k}{2}.\infty-\sum_{b\in {\mathbb{F}} }\frac{k}{2}.b&k\mbox{ even}\\\frac{k+1}{2}.\infty-\sum_{b\in {\mathbb{F}} }\frac{k-1}{2}.b&k\mbox{ odd}\end{array}\right.$$then we have a natural identification $${\mathcal L}(D)=\iota^{-1}({\cal O}_{\widetilde{\mathfrak{X}}}(-k)\otimes_{{\mathcal O}_{\widetilde{{\mathfrak X}}}}{\mathcal O}_{Z_{\gamma_0}})^c.$$In this way we may view $\iota^{-1}({\cal O}_{\widetilde{\mathfrak{X}}}(-k)\otimes_{{\mathcal O}_{\widetilde{{\mathfrak X}}}}{\mathcal O}_{Z_{\gamma_0}})^c$ as a subsheaf of ${\mathcal L}(k.\infty)$. On the other hand we may view $\iota^{-1}({\cal H}^0({\mathcal R}_{\widehat{{\mathfrak X}}}^{\bullet}))_1^{Z_{\gamma_0}}$ as a subsheaf of the constant ${\mathbb{F}} $-vector space sheaf ${\mathcal H}$ on $\mathbb{P}^1_{\mathbb{F}} $ with value $\oplus_{i=0}^k{\mathbb{F}} .z^i$ (as a sub ${\mathbb{F}} $-vector space of the function field ${\mathbb{F}} (z)$). The inclusion $\beta:{\mathcal H}\to {\mathcal L}(k.\infty)$ induces our map $\iota^{-1}({\cal H}^0({\mathcal R}_{\widehat{{\mathfrak X}}}^{\bullet}))_1^{Z_{\gamma_0}}\to{\mathcal L}(D)$ in question. It also induces an isomorphism between the respective cokernel (skyscraper) sheaves$$\frac{{\mathcal H}}{\iota^{-1}({\cal H}^0({\mathcal R}_{\widehat{{\mathfrak X}}}^{\bullet}))_1^{Z_{\gamma_0}}}\cong\frac{{\mathcal L}(k.\infty)}{{\mathcal L}(D)}$$(use the above local description of $\iota^{-1}({\cal H}^0({\mathcal R}_{\widehat{{\mathfrak X}}}^{\bullet}))_1^{Z_{\gamma_0}}$). Since clearly $\beta$ induces isomorphisms$$H^*(\mathbb{P}^1_{\mathbb{F}} ,{\mathcal H})\cong H^*(\mathbb{P}^1_{\mathbb{F}} ,{\mathcal L}(k.\infty))$$we are done.\hfill$\Box$

\begin{kor} We have the Hodge decomposition\begin{gather}H^1(\widetilde{{\mathfrak X}},{\mathcal R}_{\widehat{{\mathfrak X}}}^{\bullet})\cong H^0(\widetilde{{\mathfrak X}},{\mathcal O}_{\widehat{{\mathfrak X}}}(k+2))\oplus H^1({\widetilde{\mathfrak{X}}},{\mathcal O}_{\widehat{\mathfrak{X}}}(-k))\label{hdc}.\end{gather}
\end{kor}

{\sc Proof:} Consider the canonical maps of sheaf complexes$$[{\cal H}^0({\mathcal R}_{\widehat{{\mathfrak X}}}^{\bullet})\stackrel{0}{\to}{\mathcal R}_{\widehat{{\mathfrak X}}}^{1}]\longrightarrow{\mathcal R}_{\widehat{{\mathfrak X}}}^{\bullet}$$$$[{\cal H}^0({\mathcal R}_{\widehat{{\mathfrak X}}}^{\bullet})\stackrel{0}{\to}{\mathcal R}_{\widehat{{\mathfrak X}}}^{1}]\longrightarrow[{\mathcal R}_{\widehat{{\mathfrak X}}}^{0}\stackrel{0}{\to}{\mathcal R}_{\widehat{{\mathfrak X}}}^{1}]$$on $\widehat{{\mathfrak X}}$. By \ref{symdec} both of them induce isomorphisms in cohomology; together we thus obtain the isomorphism$$\mathbb{R}\mathbf{\Gamma}(\widetilde{{\mathfrak X}},{\mathcal R}_{\widehat{{\mathfrak X}}}^{\bullet})\cong\mathbb{R}\mathbf{\Gamma}(\widetilde{{\mathfrak X}},[{\mathcal R}_{\widehat{{\mathfrak X}}}^{0}\stackrel{0}{\to}{\mathcal R}_{\widehat{{\mathfrak X}}}^{1}]).$$We derive the stated Hodge decomposition. \hfill$\Box$\\

Let again $\Gamma<{\rm SL}\sb 2(K)$ be a cocompact discrete torsion free subgroup.

\begin{satz}\label{hodge} (a) The reduced Hodge spectral sequence$$E_1^{r,s}=H^s(X_{\Gamma},({\mathcal R}^r_X)^{\Gamma})\Rightarrow H^{r+s}(X_{\Gamma},({\mathcal R}_X^{\bullet})^{\Gamma})=H^{r+s}(X_{\Gamma},(\Omega_X^{\bullet}\otimes_K\sym_K^k({\rm St}))^{\Gamma})$$degenerates in $E_1$.\\(b) $H^{1}(X_{\Gamma},(\Omega_X^{\bullet}\otimes_K\sym_K^k({\rm St}))^{\Gamma})=H^1(X_{\Gamma},({\mathcal R}_{X}^{\bullet})^{{\Gamma}})$ decomposes naturally as$$H^{1}(X_{\Gamma},(\Omega_X^{\bullet}\otimes_K\sym_K^k({\rm St}))^{\Gamma})=H^1(\Gamma,\sym_K^k({\rm St}))\oplus H^0(X_{\Gamma},{\mathcal O}_{X}(k+2)^{\Gamma}).$$(c) If $\Gamma$ is the free group on $g$ generators and if $k>0$, then $$\dim_K(H^1(X_{\Gamma},{\cal O}_X(-k)^{\Gamma}))=\dim_K(H^0(X_{\Gamma},{\cal O}_X(k+2)^{\Gamma}))=(g-1)(k+1).$$\end{satz}

{\sc Proof:} We may of course change bases from $K$ to $\widehat{K}$. Statement (a) is a consequence of \ref{gammakoh} (if $k>0$) but we can also argue as follows. It is enough to show that the inclusion of sheaf complexes $$[{\cal H}^0({\mathcal R}^{\bullet}_{\widehat{X}})^{\Gamma}\stackrel{0}{\longrightarrow}({\mathcal R}_{\widehat{X}}^{1})^{\Gamma}]\hookrightarrow({\mathcal R}^{\bullet}_{\widehat{X}})^{\Gamma}$$on $\widehat{X}_{\Gamma}$ induces isomorphisms $$H^*(\widetilde{X}_{\Gamma},[{\cal H}^0({\mathcal R}^{\bullet}_{\widehat{X}})^{\Gamma}\stackrel{0}{\longrightarrow}({\mathcal R}_{\widehat{X}}^{1})^{\Gamma}])\cong H^*(\widetilde{X}_{\Gamma},({\mathcal R}^{\bullet}_{\widehat{X}})^{\Gamma}).$$For this it suffices to show that ${\cal H}^0({\mathcal R}^{\bullet}_{\widehat{X}})^{\Gamma}\to({\mathcal R}_{\widehat{X}}^{0})^{\Gamma}$ induces isomorphisms in cohomology. Now $\widehat{X}_{\Gamma}$ is quasi-compact, hence $H^*(\widetilde{X}_{\Gamma},.)$ commutes with $(.)\otimes_{{\cal O}_{\widehat{K}}}\widehat{K}$. Therefore it suffices to show that the morphism of sheaves ${\cal H}^0({\mathcal R}^{\bullet}_{\widehat{\mathfrak{X}}})^{\Gamma}\to({\mathcal R}_{\widehat{\mathfrak{X}}}^{0})^{\Gamma}$ on $\widehat{\mathfrak{X}}_{\Gamma}$ induces isomorphisms \begin{gather}H^*(\widetilde{\mathfrak{X}}_{\Gamma},{\cal H}^0({\mathcal R}^{\bullet}_{\widehat{\mathfrak{X}}})^{\Gamma})\cong H^*(\widetilde{\mathfrak{X}}_{\Gamma},({\mathcal R}_{\widehat{\mathfrak{X}}}^{0})^{\Gamma})\label{vvvvv}.\end{gather} Using the covering spectral sequences$$E_2^{r,s}=H^r(\Gamma,H^s(\widehat{\mathfrak{X}},{\cal F}))\Rightarrow H^{r+s}(\widehat{\mathfrak{X}}_{\Gamma},{\cal F}^{\Gamma})$$for ${\cal F}={\cal H}^0({\mathcal R}^{\bullet}_{\widehat{\mathfrak{X}}})$ and ${\cal F}={\mathcal R}_{\widehat{\mathfrak{X}}}^{0}$ we see that it is enough to prove that the maps$$H^r(\Gamma,H^s({\widetilde{\mathfrak{X}}},{\cal H}^0({\mathcal R}^{\bullet}_{\widehat{\mathfrak{X}}})))\to H^r(\Gamma,H^s({\widetilde{\mathfrak{X}}},{\mathcal R}_{\widehat{\mathfrak{X}}}^{0}))$$are isomorphisms. But they are, as follows from \ref{symdec}. We turn to (b). We have \begin{align}H^1(\widetilde{X}_{\Gamma},({\mathcal R}_{\widehat{X}}^{\bullet})^{{\Gamma}})&=H^1(\widetilde{\mathfrak{X}}_{\Gamma},({\mathcal R}^{\bullet}_{\widehat{\mathfrak{X}}})^{\Gamma})\otimes_{{\cal O}_{\widehat{K}}}\widehat{K}\notag \\
{} & =H^1(\widetilde{\mathfrak{X}}_{\Gamma},{\cal H}^0({\mathcal R}^{\bullet}_{\widehat{\mathfrak{X}}})^{\Gamma})\otimes_{{\cal O}_{\widehat{K}}}\widehat{K}\oplus H^0(\widetilde{\mathfrak{X}}_{\Gamma},({\mathcal R}^{1}_{\widehat{\mathfrak{X}}})^{\Gamma})\otimes_{{\cal O}_{\widehat{K}}}\widehat{K}\notag \\
{} & =H^1(\widehat{{X}}_{\Gamma},{\cal H}^0({\mathcal R}^{\bullet}_{\widehat{{X}}})^{\Gamma})\oplus H^0(\widehat{{X}}_{\Gamma},({\mathcal R}^{1}_{\widehat{{X}}})^{\Gamma})\notag\end{align}where the first and the third equality follow again from the quasi-compactness of $\widehat{X}_{\Gamma}$, and the second equality from (\ref{vvvvv}). Now ${\mathcal R}^{1}_{{{X}}}={\mathcal O}_{X}(k+2)$, and on the other hand$$H^1({{X}}_{\Gamma},{\cal H}^0({\mathcal R}^{\bullet}_{{{X}}})^{\Gamma})=H^1(\Gamma,\mathbb{R}\mathbf{\Gamma}(X,{\cal H}^0({\mathcal R}^{\bullet}_{{{X}}}))).$$But $H^0(X,{\cal H}^0({\mathcal R}^{\bullet}_{{{X}}}))=H^0(X,{\mathcal R}^{\bullet}_{{{X}}})=\sym_K^k({\rm St})$ and $H^j(X,{\cal H}^0({\mathcal R}^{\bullet}_{{{X}}}))=0$ for $j\ne0$ because ${\cal H}^0({\mathcal R}^{\bullet}_{{{X}}})$ is the locally constant sheaf on $X$ generated by $H^0(X,{\mathcal R}^{\bullet}_{{{X}}})=H^0(X,\Omega_X^{\bullet}\otimes_K\sym_K^k({\rm St}))=\sym_K^k({\rm St})$. In (c) for the equality $\dim_K(H^0(X_{\Gamma},{\cal O}_X(k+2)^{\Gamma}))=(g-1)(k+1)$ see \cite{ss} p.98. The equality $\dim_K(H^1(X_{\Gamma},{\cal O}_X(-k)^{\Gamma}))=(g-1)(k+1)$ follows from statement (a) together with \cite{schn} p.628 and \cite{ss} p.98.\hfill$\Box$

The decomposition in (b) is not new. It was established for the first time in \cite{deshhod} and later again in \cite{schn}. Both these (mutually different) proofs use sophisticated analytic methods (e.g. Coleman integration in \cite{deshhod}). The degeneration of the spectral sequence in (a) however, conjectured in \cite{schn}, seemed to be unknown before (cf. \cite{schn} p.649). Note that the spectral sequence for the {\it non-reduced} de Rham complex does not degenerate in general at $E_1$ (cf. loc. cit.). 

\begin{kor} The intersection of $H^0(\widetilde{\mathfrak X},{\mathcal O}_{\widehat{\mathfrak X}}(k+2))$ and of$$\bi[H^0(\widetilde{X},{\mathcal O}_{\widehat{X}}(-k))\stackrel{\partial^{k+1}}{\longrightarrow}H^0(\widetilde{X},{\mathcal O}_{\widehat{X}}(k+2))]$$inside $H^0(\widetilde{\mathfrak X},sp_*{\mathcal O}_{\widehat{X}}(k+2))=H^0(\widetilde{X},{\mathcal O}_{\widehat{X}}(k+2 ))$ is zero. In particular, $H^0(\widetilde{\mathfrak X},{\mathcal O}_{\widehat{\mathfrak X}}(k+2))$ can be viewed as a submodule of $H^1(\widetilde{X},{\mathcal R}_{\widehat{X}}^{\bullet})=H^1(\widetilde{X},\Omega_{\widehat{X}}^{\bullet}\otimes_K\sym_K^k({\rm St}))$.
\end{kor}

{\sc Proof:} This follows immediately from the injectivity of the map $Res^0$ in \ref{harm}.\hfill$\Box$ 

\begin{satz}\label{cechflip} For $k>0$ there is a natural $G$-equivariant isomorphism$$\theta:Z_{har}^1(k+2)\otimes\varepsilon^{k+1}\cong H^1({\widetilde{\mathfrak{X}}},{\cal H}^0({\mathcal R}^{\bullet}_{\widehat{\mathfrak{X}}})).$$
\end{satz}

{\sc Proof:} Observe $\varepsilon=\det\cdot\chi^{-2}$, which relates the twisting here to that in Section \ref{harsec}. We have a $G$-equivariant isomorphism$$\Hom_{\widehat{K}}(\sym^k_{\widehat{K}}({\rm St})[-k]\otimes\chi^k,\widehat{K})\stackrel{\sigma}{\longrightarrow}H^0(\widetilde{X},{\cal H}^0({\mathcal R}^{\bullet}_{\widehat{{X}}})),\quad h_{j}\mapsto z^{k-j}$$with $h_j\in\Hom_{\widehat{K}}(\sym^k_{\widehat{K}}({\rm St})[-k]\otimes\chi^k,{\widehat{K}})$ as in the proof of \ref{harm}, i.e. $h_j(X^jY^{k-j})=1$ and $h_j(X^iY^{k-i})=0$ for $i\ne j$. For $Z\in F^0$ and $\{Z_1,Z_2\}\in F^1$ we define sheaves ${\cal G}_Z$ and ${\cal G}_{\{Z_1,Z_2\}}$ on ${\widehat{\mathfrak{X}}}$: for open $U\subset{\widehat{\mathfrak{X}}}$ we let$${\cal G}_Z(U)=\left\{\begin{array}{l@{\quad:\quad}l}{\cal H}^0({\mathcal R}^{\bullet}_{\widehat{\mathfrak{X}}})(\widehat{\mathfrak U}_{\{Z\}})&U\cap Z\ne\emptyset\\0&U\cap Z=\emptyset\end{array}\right.$$$${\cal G}_{\{Z_1,Z_2\}}(U)=\left\{\begin{array}{l@{\quad:\quad}l}{\cal G}_{Z_1}(U)+{\cal G}_{Z_2}(U)&U\cap Z_1\cap Z_2\ne\emptyset\\0&U\cap Z_1\cap Z_2=\emptyset\end{array}\right..$$Then we have an exact sequence\begin{gather}0\longrightarrow{\cal H}^0({\mathcal R}^{\bullet}_{\widehat{\mathfrak{X}}})\longrightarrow\prod_{Z\in F^0}{\cal G}_Z\stackrel{\delta}{\longrightarrow}\prod_{\{Z_1,Z_2\}\in F^1}{\cal G}_{\{Z_1,Z_2\}}\longrightarrow0\label{ceha}\end{gather}where $\delta$ is the product of all maps ${\rm sg}(Z_1).\id:{\cal G}_{Z_1}\to{\cal G}_{\{Z_1,Z_2\}}$. In cohomology we get$$\frac{H^0({\widetilde{\mathfrak{X}}},\prod_{\{Z_1,Z_2\}\in F^1}{\cal G}_{\{Z_1,Z_2\}})}{H^0({\widetilde{\mathfrak{X}}},\prod_{Z\in F^0}{\cal G}_Z)}\cong H^1({\widetilde{\mathfrak{X}}},{\cal H}^0({\mathcal R}^{\bullet}_{\widehat{\mathfrak{X}}})).$$We claim that$$Z_{har}^1(k+2)\otimes\varepsilon^{k+1}\to H^0({\widetilde{\mathfrak{X}}},\prod_{\{Z_1,Z_2\}\in F^1}{\cal G}_{\{Z_1,Z_2\}})$$$$(f_{\{Z_1,Z_2\}})_{\{Z_1,Z_2\}}\mapsto\prod_{\{Z_1,Z_2\}}\sigma(f_{\{Z_1,Z_2\}})$$induces an isomorphism $\theta:Z_{har}^1(k+2)\otimes\varepsilon^{k+1}\to H^1({\widetilde{\mathfrak{X}}},{\cal H}^0({\mathcal R}^{\bullet}_{\widehat{\mathfrak{X}}}))$. Since $H^1({\widetilde{\mathfrak{X}}},{\cal H}^0({\mathcal R}^{\bullet}_{\widehat{\mathfrak{X}}}))=H^1({\widetilde{\mathfrak{X}}},{\mathcal O}_{\widehat{{\mathfrak X}}}(-k))$ is flat it suffices to show that the induced map $$\widetilde{\theta}=\theta/(\widehat{\pi}):\widetilde{Z}_{har}^1(k+2)\longrightarrow H^1({\widetilde{\mathfrak{X}}},({\cal H}^0({\mathcal R}^{\bullet}_{\widehat{\mathfrak{X}}}))_1)$$is an isomorphism (with notations from the proof of \ref{symdec}). Let us first assume $k>0$ is even. Consider the submodule$$\widetilde{Z}_{har}^1(k+2)(1)=\left\{\begin{array}{l@{\quad}l}f=(f_{\{Z_1,Z_2\}})_{\{Z_1,Z_2\}\in F^1}\in\widetilde{Z}_{har}^1(k+2);\\(\gamma.f)_{\{Z_{\gamma_0},Z_{\gamma_1}\}}(X^{\frac{k}{2}}Y^{\frac{k}{2}})=0\mbox{ for all }\gamma\in G\end{array}\right\}$$of $\widetilde{Z}_{har}^1(k+2)$ (this is nothing but the image of $H^0(\widetilde{\mathfrak{X}},{\mathcal O}_{\widetilde{{\mathfrak X}}}(k+2)(1))$ under $\widetilde{Res}^0$). If for $Z\in F^0$ we let $\widetilde{Z}_{har}^1(k+2)(1)_Z$ be the image of $\widetilde{Z}_{har}^1(k+2)(1)\to\widetilde{Z}_{har}^1(k+2)\to\widetilde{Z}_{har}^1(k+2)_Z$, then$$\widetilde{Z}_{har}^1(k+2)(1)=\prod_{Z\in F^0}\widetilde{Z}_{har}^1(k+2)(1)_Z.$$ In particular we have natural injections $\iota_Z:\widetilde{Z}_{har}^1(k+2)(1)_Z\to\widetilde{Z}_{har}^1(k+2)(1)$. We claim that for each $Z\in F^0$ the composition$$\widetilde{Z}_{har}^1(k+2)(1)_Z\stackrel{\iota_Z}{\longrightarrow}\widetilde{Z}_{har}^1(k+2)(1)\stackrel{\widetilde{\theta}}{\longrightarrow}H^1({\widetilde{\mathfrak{X}}},({\cal H}^0({\mathcal R}^{\bullet}_{\widehat{\mathfrak{X}}}))_1)\longrightarrow H^1({\widetilde{\mathfrak{X}}},({\cal H}^0({\mathcal R}^{\bullet}_{\widehat{\mathfrak{X}}}))^Z_1),$$which we denote by $\beta_Z$, is an isomorphism. To see this we may assume $Z=Z_{\gamma_0}$. From the proof of \ref{symdec} we infer an exact sequence$$0\longrightarrow H^0(\mathbb{P}^1,{\cal H})\stackrel{}{\longrightarrow}H^0(\mathbb{P}^1,\frac{{\cal H}}{\iota^{-1}({\cal H}^0({\mathcal R}^{\bullet}_{\widehat{\mathfrak{X}}}))^{Z_{\gamma_0}}_1})\longrightarrow H^1({\widetilde{\mathfrak{X}}},({\cal H}^0({\mathcal R}^{\bullet}_{\widehat{\mathfrak{X}}}))_1^{Z_{\gamma_0}})\longrightarrow0.$$Here $\iota:\mathbb{P}^1\cong Z_{\gamma_0}\to{\widetilde{\mathfrak{X}}}$ is the natural embedding, ${\cal H}$ is the constant sheaf with value $\oplus_{i=0}^k{\mathbb{F}} .z^i$, and the quotient ${\cal H}/\iota^{-1}({\cal H}^0({\mathcal R}^{\bullet}_{\widehat{\mathfrak{X}}}))^{Z_{\gamma_0}}_1$ is a skyscraper sheaf whose only stalks are $\frac{k}{2}$-dimensional ${\mathbb{F}} $-vector spaces at the ${\mathbb{F}} $-rational points of $\mathbb{P}^1$. Namely, in notations from section \ref{cosec}, the ${\mathbb{F}} $-rational points of $\mathbb{P}^1\cong Z_{\gamma_0}$ are just the intersections $Z_{\gamma_0}\cap Z_{\gamma_{a,0}\gamma_{-1}}$ with $a\in R$, and $Z_{\gamma_0}\cap Z_{\gamma_{1}}$. The stalk of ${\cal H}/\iota^{-1}({\cal H}^0({\mathcal R}^{\bullet}_{\widehat{\mathfrak{X}}}))^{Z_{\gamma_0}}_1$ at $Z_{\gamma_0}\cap Z_{\gamma_{a,0}\gamma_{-1}}$ is (canonically identified with) $\oplus_{i=0}^{\frac{k}{2}-1}{\mathbb{F}} .(z-\overline{a})^{i}$ (with $\overline{a}\in {\mathbb{F}} $ the image of $a\in R$), and the stalk at $Z_{\gamma_0}\cap Z_{\gamma_{1}}$ is (canonically identified with) $\oplus_{i=\frac{k}{2}+1}^{k}{\mathbb{F}} .z^{i}$. From the proof of \ref{harm} we get the exact sequence$$0\longrightarrow\widetilde{Z}_{har}^1(k+2)(1)_{Z_{\gamma_0}}\longrightarrow\oplus_{j=0}^{\frac{k}{2}-1}{\mathbb{F}} .h_j\times\prod_{a\in R}\oplus_{j=\frac{k}{2}+1}^{k}{\mathbb{F}} .(\gamma_{a,0}h_{j})\stackrel{\sum}{\longrightarrow}\oplus_{j=0}^{{k}}{\mathbb{F}} .h_j$$(the first factor in the middle term is the  $\{Z_{\gamma_0},Z_{\gamma_{1}}\}$-component). Now $\sigma$ maps $\gamma_{a,0}.h_j$ to $\gamma_{a,0}.z^{k-j}=(z-\overline{a})^{k-j}$, hence defines a map$$\widetilde{Z}_{har}^1(k+2)(1)_{Z_{\gamma_0}}\to H^0(\mathbb{P}^1,\frac{{\cal H}}{\iota^{-1}({\cal H}^0({\mathcal R}^{\bullet}_{\widehat{\mathfrak{X}}}))^{Z_{\gamma_0}}_1})$$whose composition with the projection to $H^1({\widetilde{\mathfrak{X}}},({\cal H}^0({\mathcal R}^{\bullet}_{\widehat{\mathfrak{X}}}))^{Z_{\gamma_0}}_1)$ is an isomorphism: this isomorphism is our $\beta_{Z_{\gamma_0}}$. We have shown that $\widetilde{\theta}|_{\widetilde{Z}_{har}^1(k+2)(1)}$ is injective and that its image $\bi(\widetilde{\theta}|_{\widetilde{Z}_{har}^1(k+2)(1)})$$\subset H^1({\widetilde{\mathfrak{X}}},({\cal H}^0({\mathcal R}^{\bullet}_{\widehat{\mathfrak{X}}}))_1)$ maps isomorphically to $H^1({\widetilde{\mathfrak{X}}},\prod_{Z\in F^0}({\cal H}^0({\mathcal R}^{\bullet}_{\widehat{\mathfrak{X}}}))^Z_1)$. From the exact sequence$$0\to H^0(\widetilde{\mathfrak{X}},\prod_{\{Z_1,Z_2\}\in F^1}({\cal H}^0({\mathcal R}^{\bullet}_{\widehat{\mathfrak{X}}}))_1^{Z_1,Z_2})\longrightarrow H^1({\widetilde{\mathfrak{X}}},({\cal H}^0({\mathcal R}^{\bullet}_{\widehat{\mathfrak{X}}}))_1)\longrightarrow H^1({\widetilde{\mathfrak{X}}},\prod_{Z\in F^0}({\cal H}^0({\mathcal R}^{\bullet}_{\widehat{\mathfrak{X}}}))^Z_1)\longrightarrow0$$we therefore get$$\frac{H^1({\widetilde{\mathfrak{X}}},({\cal H}^0({\mathcal R}^{\bullet}_{\widehat{\mathfrak{X}}}))_1)}{\bi(\widetilde{\theta}|_{\widetilde{Z}_{har}^1(k+2)(1)})}\cong H^0(\widetilde{\mathfrak{X}},\prod_{\{Z_1,Z_2\}\in F^1}({\cal H}^0({\mathcal R}^{\bullet}_{\widehat{\mathfrak{X}}}))_1^{Z_1,Z_2}).$$In particular we get a map$$\frac{\widetilde{Z}_{har}^1(k+2)}{\widetilde{Z}_{har}^1(k+2)(1)}\longrightarrow H^0(\widetilde{\mathfrak{X}},\prod_{\{Z_1,Z_2\}\in F^1}({\cal H}^0({\mathcal R}^{\bullet}_{\widehat{\mathfrak{X}}}))_1^{Z_1,Z_2})$$induced by $\widetilde{\theta}$ and it remains to show that this map is bijective. But this is clear, as both sides can be identified with (\ref{x}). If $k>0$ is odd things are easier since there are no terms $({\cal H}^0({\mathcal R}^{\bullet}_{\widehat{\mathfrak{X}}}))_1^{Z_1,Z_2}$ and we only need to show bijectivity of the maps$$\widetilde{Z}_{har}^1(k+2)_Z\stackrel{\iota_Z}{\longrightarrow}\widetilde{Z}_{har}^1(k+2)\stackrel{\widetilde{\theta}}{\longrightarrow}H^1({\widetilde{\mathfrak{X}}},({\cal H}^0({\mathcal R}^{\bullet}_{\widehat{\mathfrak{X}}}))_1)\longrightarrow H^1({\widetilde{\mathfrak{X}}},({\cal H}^0({\mathcal R}^{\bullet}_{\widehat{\mathfrak{X}}}))^Z_1).$$We can proceed just as before, now the sums in our local analysis run from $0$ to $\frac{k-1}{2}$, resp. from $\frac{k+1}{2}$ to $k$.\hfill$\Box$\\

At this point we see that by considering integral structures in our automorphic line bundles ${\mathcal{O}}_X(k)$ on $X$ we obtain genuinely new structures in cohomology. Namely, whereas Theorem \ref{harm} does have a non-integral counterpart --- the isomorphism$$Res:\frac{\Gamma(X,{\mathcal O}_X(k+2))}{\bi[\Gamma(X,{\mathcal O}_X(-k))\stackrel{{\partial}^{k+1}}{\longrightarrow}\Gamma(X,{\mathcal O}_X(k+2))]}\cong C^1_{har}(K)$$from \cite{ss} p.97 ---, Theorem \ref{cechflip} has no non-integral counterpart (in fact $H^1(X,{\cal H}^0({\mathcal R}^{\bullet}_X))=0$). As an application of Theorem \ref{cechflip} we get a {\it global} version of the monodromy operator, as follows. From \ref{harm}, \ref{symdec} and \ref{cechflip} we obtain $G$-equivariant isomorphisms (if $k>0$)$$H^0(\widetilde{{\mathfrak X}},{\mathcal O}_{\widehat{{\mathfrak X}}}(k+2))\cong Z_{har}^1(k+2)\cong H^1({\widetilde{\mathfrak{X}}},{\cal H}^0({\mathcal R}^{\bullet}_{\widehat{\mathfrak{X}}}))\otimes\varepsilon^{-k-1}\cong  H^1({\widetilde{\mathfrak{X}}},{\mathcal O}_{\widehat{\mathfrak{X}}}(-k))\otimes\varepsilon^{-k-1}$$whose composition we denote by $\nu$.

{\bf Definition:} The monodromy operator $N:H^1(\widetilde{{\mathfrak X}},{\mathcal R}_{\widehat{{\mathfrak X}}}^{\bullet})\to H^1(\widetilde{{\mathfrak X}},{\mathcal R}_{\widehat{{\mathfrak X}}}^{\bullet})$ is the composition$$H^1(\widetilde{{\mathfrak X}},{\mathcal R}_{\widehat{{\mathfrak X}}}^{\bullet})\stackrel{pr}{\longrightarrow}H^0(\widetilde{{\mathfrak X}},{\mathcal O}_{\widehat{{\mathfrak X}}}(k+2))\stackrel{\nu}{\longrightarrow}H^1({\widetilde{\mathfrak{X}}},{\mathcal O}_{\widehat{\mathfrak{X}}}(-k))\stackrel{i}{\longrightarrow}H^1(\widetilde{{\mathfrak X}},{\mathcal R}_{\widehat{{\mathfrak X}}}^{\bullet})$$where $pr$ resp. $i$ is the natural projection resp. inclusion in (\ref{hdc}).

Thus $N$ is $G$-equivariant when viewed as a map $H^1(\widetilde{{\mathfrak X}},{\mathcal R}_{\widehat{{\mathfrak X}}}^{\bullet})\to H^1(\widetilde{{\mathfrak X}},{\mathcal R}_{\widehat{{\mathfrak X}}}^{\bullet})\otimes\varepsilon^{-k-1}$. Its monodromy filtration $\ke(N)=\bi(N)=H^1({\widetilde{\mathfrak{X}}},{\mathcal O}_{\widehat{\mathfrak{X}}}(-k))$ splits the Hodge filtration $H^0(\widetilde{{\mathfrak X}},{\mathcal O}_{\widehat{{\mathfrak X}}}(k+2))$ of $H^1(\widetilde{{\mathfrak X}},{\mathcal R}_{\widehat{{\mathfrak X}}}^{\bullet})$. Now we restrict our attention to the action by ${\rm SL}\sb 2(K)$. If $\Gamma<{\rm SL}\sb 2(K)$ is a cocompact discrete torsion free subgroup, we only need to take $\Gamma$-invariants and invert $p$ in (\ref{hdc}) to obtain the Hodge decomposition$$H^{1}(\widehat{X}_{\Gamma},(\Omega_{\widehat{X}}^{\bullet}\otimes_{\widehat{K}}\sym_{\widehat{K}}^k({\rm St}))^{\Gamma})= H^0(\widetilde{X}_{\Gamma},{\mathcal O}_{{{\widehat{X}}}}(k+2)^{\Gamma})\oplus H^1(\Gamma,\sym_{\widehat{K}}^k({\rm St}))$$from \ref{hodge} (we saw $H^1(\Gamma,\sym_{\widehat{K}}^k({\rm St}))=H^1({\widetilde{\mathfrak{X}}},{\mathcal O}_{\widehat{\mathfrak{X}}}(-k))^{\Gamma}\otimes{\mathbb{Q}}$ in \ref{hodge}): no higher $\Gamma$-group cohomology is needed for this passage. It is not hard to see that the monodromy operator we thus obtain on $H^{1}(\widehat{X}_{\Gamma},(\Omega_{\widehat{X}}^{\bullet}\otimes_{\widehat{K}}\sym_{\widehat{K}}^k({\rm St}))^{\Gamma})$ is the one predicted by $p$-adic Hodge theory, using the description of the latter given in \cite{iovspi}. In particular this shows that $N$ respects the integral {\it de Rham} structures (as opposed to integral Hyodo-Kato cohomology structures) in $H^{1}(\widehat{X}_{\Gamma},(\Omega_{\widehat{X}}^{\bullet}\otimes_{\widehat{K}}\sym_{\widehat{K}}^k({\rm St}))^{\Gamma})$, a fact which the general $p$-adic Hodge theory does not seem to suggest. We so obtain an infinite rank filtered monodromy module over ${\cal O}_{\widehat{K}}$ which comprises all the filtered monodromy modules $H^{1}(\widehat{X}_{\Gamma},(\Omega_{\widehat{X}}^{\bullet}\otimes_{\widehat{K}}\sym_{\widehat{K}}^k({\rm St}))^{\Gamma})$ for the various $\Gamma$.

For $k=0$ we still can define $N$ de Rham integrally as the composition$$H^1(\widetilde{{\mathfrak X}}_{\Gamma},({\mathcal R}_{\widehat{{\mathfrak X}}}^{\bullet})^{\Gamma})\stackrel{pr}{\longrightarrow}H^0(\widetilde{{\mathfrak X}}_{\Gamma},{\mathcal O}_{\widehat{{\mathfrak X}}}(2)^{\Gamma})\stackrel{Res^0}{\longrightarrow}Z_{har}^1(2)^{\Gamma}$$$$\stackrel{\xi}{\longrightarrow} H^0({\widetilde{\mathfrak{X}}},\prod_{\{Z_1,Z_2\}}{\cal G}_{\{Z_1,Z_2\}})^{\Gamma}\stackrel{\delta}{\longrightarrow} H^1(\Gamma,H^0({\widetilde{\mathfrak{X}}},{\cal H}^0({\mathcal R}^{\bullet}_{\widehat{\mathfrak{X}}})))\stackrel{i}{\longrightarrow}H^1(\widetilde{{\mathfrak X}}_{\Gamma},({\mathcal R}_{\widehat{{\mathfrak X}}}^{\bullet})^{\Gamma}).$$Here sheaves ${\cal G}_Z$ and ${\cal G}_{\{Z_1,Z_2\}}$ and a map $Z_{har}^1(2)\to H^0({\widetilde{\mathfrak{X}}},\prod_{\{Z_1,Z_2\}}{\cal G}_{\{Z_1,Z_2\}})$ are defined just as in the proof of \ref{cechflip}, and $\xi$ is the restricted map on $\Gamma$-invariants. The map $\delta$ is the connecting homomorphism in group cohomology (observe that for $k=0$ application of $H^0(\widetilde{\mathfrak{X}},.)$ to the sequence (\ref{ceha}) preserves its exactness). Inverting $p$ in the above composition gives the correct $N$ on $H^{1}(\widehat{X}_{\Gamma},\Omega_{\widehat{X}_{\Gamma}}^{\bullet})$ (at least up to sign, see \cite{iovspi}).

\section{Complements}

(A) Let ${\mathcal R}_{\widetilde{\mathfrak X}}^{\bullet}={\mathcal R}_{\widehat{\mathfrak X}}^{\bullet}/(\widehat{\pi})$. One can prove the analogs of \ref{symdec} and \ref{hodge} for ${\mathcal R}_{\widetilde{\mathfrak X}}^{\bullet}$, namely:\\
$$H^1(\widetilde{\mathfrak X},{\mathcal R}_{\widetilde{\mathfrak X}}^{\bullet})=H^1(\widetilde{\mathfrak X},{\cal H}^0({\mathcal R}_{\widetilde{\mathfrak X}}^{\bullet}))\oplus H^0(\widetilde{\mathfrak X},{\mathcal O}_{\widetilde{\mathfrak X}}(k+2))$$$$H^1(\widetilde{\mathfrak X}_{\Gamma},({\mathcal R}_{\widetilde{\mathfrak X}}^{\bullet})^{\Gamma})=H^1(\widetilde{\mathfrak X}_{\Gamma},({\cal H}^0({\mathcal R}_{\widetilde{\mathfrak X}}^{\bullet}))^{\Gamma})\oplus H^0(\widetilde{\mathfrak X}_{\Gamma},{\mathcal O}_{\widetilde{\mathfrak X}}(k+2)^{\Gamma}).$$Note that this is {\it not} obvious from the proof of \ref{symdec}, there we did {\it not} consider ${\cal H}^0({\mathcal R}_{\widetilde{\mathfrak X}}^{\bullet})$.\\ 

(B) Let $k\in\mathbb{Z}$ be even and let $\omega_{\widehat{\mathfrak X}/{\mathcal O}_{\widehat{K}}}$ be the logarithmic differential module of the log smooth morphism $\widehat{\mathfrak X}\to\spf({\mathcal O}_{\widehat{K}})$: an invertible ${\rm PGL}\sb 2(K)$-equivariant line bundle on $\widehat{\mathfrak X}$. We have an ${\rm SL}\sb 2(K)$-equivariant isomorphism$${\mathcal O}_{\widehat{\mathfrak X}}(k)\cong\omega_{\widehat{\mathfrak X}/{\mathcal O}_{\widehat{K}}}^{\frac{k}{2}},\quad f\mapsto fdz^{\frac{k}{2}}.$$Now $dz^{\frac{k}{2}}$ is not a generator of $\omega^{\frac{k}{2}}_{\widehat{\mathfrak X}/{\mathcal O}_{\widehat{K}}}$, not even a global section of $\omega^{\frac{k}{2}}_{\widehat{\mathfrak X}/{\mathcal O}_{\widehat{K}}}$ if $k<0$. Let $k>0$ and even. For $a\in{\mathcal O}_K$ the local section $\dlog(z-a)$ is a generator of $\omega_{\widehat{\mathfrak X}/{\mathcal O}_{\widehat{K}}}$ on an appropriate open formal subscheme of $\widehat{\mathfrak X}$. There, the complex ${\mathcal R}_{\widehat{\mathfrak X}}^{\bullet}$ becomes isomorphic to$$\omega_{\widehat{\mathfrak X}/{\mathcal O}_{\widehat{K}}}^{\frac{-k}{2}}\longrightarrow\omega_{\widehat{\mathfrak X}/{\mathcal O}_{\widehat{K}}}^{\frac{k+2}{2}}$$$$f\dlog(z-a)^{\frac{-k}{2}}\mapsto(D_a\prod_{j=1}^{\frac{k}{2}}(D_a^2-j^2)f)\dlog(z-a)^{\frac{k+2}{2}}$$where $D_a=(z-a)\partial=\frac{(z-a)d}{d(z-a)}$. For the proof you need to show $(z-a)^{\frac{k+2}{2}}\partial^{k+1}(z-a)^{\frac{k}{2}}=D_a\prod_{j=1}^{\frac{k}{2}}(D_a^2-j^2)$. For this show by induction on $n$, departing from $D_a=\partial(z-a)-1$ that $(z-a)^n\partial^n=D_a(D_a-1)\ldots(D_a-n+1)$ and $\partial^n(z-a)^n=(D_a+n)(D_a+n-1)\ldots(D_a+1)$. Also note $-D_a=(z-a)^{-1}\frac{d}{d(z-a)^{-1}}$.\\

(C) For even weights $k\in\mathbb{Z}$ the ${\cal O}_{\widehat{\mathfrak{X}}}$-modules ${\cal O}_{\widehat{\mathfrak{X}}}(k)$ are in fact line bundles, and the base extension $K\to\widehat{K}$ is unnecessary, i.e. everything we did here descends from $\widehat{\mathfrak{X}}$ to $\mathfrak{X}$. The automorphic action of even weight $k$ in \cite{jer} is the one we get by replacing the factor $\chi^k(\gamma)$ with the factor $\det(\gamma)^{\frac{k}{2}}$ in equation (\ref{xiii}). All our results carry over to this situation (and in \ref{cechflip} no $\varepsilon^{k+1}$-twist is needed). But also if the weight $k$ is odd, if one is willing to restrict the automorphic action on ${\cal O}_{{{X}}}(k)$ to a smaller group, the base extension $K\to\widehat{K}$ can be avoided and one has equivariant integral structures which are even line bundles. 
Let $G^{even}=\{\gamma\in G;\,\omega(\det(\gamma))\mbox{ even}\}$. Note that the restriction to $G^{even}$ of the automorphic action (defined in equation (\ref{xiii})) only depends on the choice of $\pi$, not of $\widehat{\pi}$. In notations from section \ref{fullgsec}, define the following ${\mathcal O}_{\mathfrak{U}_{\{Z_{\gamma_n},Z_{\gamma_{n+1}}\}}}$-submodule of ${\mathcal O}_{\mathfrak{U}_{\{Z_{\gamma_n},Z_{\gamma_{n+1}}\}}}\otimes_{{\mathcal O}_K}K$:$${\mathcal O}_{\mathfrak{U}_{\{Z_{\gamma_n},Z_{\gamma_{n+1}}\}}}(k)={\mathcal O}_{\mathfrak{U}_{\{Z_{\gamma_n},Z_{\gamma_{n+1}}\}}}.f_{n,n}^{\lfloor\frac{kn}{2}\rfloor}f_{n,n+1}^{\lfloor\frac{k(n+1)}{2}\rfloor}.$$The ${\mathcal O}_{\mathfrak{U}_{\{Z_{\gamma_n},Z_{\gamma_{n+1}}\}}}(k)$ glue into an invertible ${\mathcal O}_{\mathfrak{Y}}$-submodule ${\mathcal O}_{\mathfrak{Y}}(k)$ of ${\mathcal O}_{\mathfrak{Y}}\otimes_{{\mathcal O}_K}K$. Observe\begin{gather}{\mathcal O}_{\mathfrak{Y}}(k)|_{{\mathfrak U}_{\{Z_{\gamma_n}\}}}=\pi^{\lfloor\frac{kn}{2}\rfloor}{\mathcal O}_{{\mathfrak U}_{\{Z_{\gamma_n}\}}}\quad\mbox{inside}\quad{\mathcal O}_{{\mathfrak U}_{\{Z_{\gamma_n}\}}}\otimes_{{\mathcal O}_K}K.\notag\end{gather}As in \ref{gequivint} one sees that ${\mathcal O}_{\mathfrak{Y}}(k)$ globalizes to a $G^{even}$-equivariant line bundle ${\mathcal O}_{\mathfrak{X}}(k)$ on ${\mathfrak{X}}$, an integral structure in ${\cal O}_{{{X}}}(k)$. Our entire analysis of ${\cal O}_{\widehat{\mathfrak{X}}}(k)$ can be repeated with ${\mathcal O}_{\mathfrak{X}}(k)$, with essentially the same results (e.g. those from section \ref{cosec}; however, the case $k=1$ is slightly harder in this context). One additional feature is that one has to study ${\mathcal O}_{\mathfrak{X}}(k)\otimes_{{\mathcal O}_{\mathfrak{X}}}{\mathcal O}_{Z}$ for $Z\in F^0_{even}$ and for $Z\in F^0_{odd}$ separately (the two orbits of $G^{even}$ acting on $F^0$) and the shapes of these two types are indeed different if $k$ is odd.


\begin{flushleft}
\textsc{Mathematisches Institut der Universit\"at M\"unster\\ Einsteinstrasse 62, 48149 M\"unster, Germany}\\
\textit{E-mail address}: klonne@math.uni-muenster.de
\end{flushleft}
\end{document}